\documentclass[final]{siamltex}
\usepackage{amsfonts}
\usepackage{graphicx}
\usepackage{verbatim}

\newcommand{\C}{\mathbb{C}}

\font\ff=eusm10 scaled 1200

\def\K{\hbox{\ff K}}

\def\M{\mbox{\sf M}} 
\newcommand{\re}{\mathrm{Re}\,}

\usepackage{subfig}

\graphicspath{{PDF/}{EPS/}}

\usepackage{amsmath}

\def\IH{\mathbb{H}}
\def\ID{\mathbb{D}}

\def\Cap{\mathrm{cap}}

\usepackage{fancybox}

\newtheorem{remark}[theorem]{Remark}

\title{Conformal modulus on domains with strong singularities and cusps}

\author{Harri Hakula\thanks{Aalto University, Institute of Mathematics,
P.O. Box 11100, FI-00076 Aalto,
FINLAND ({\tt harri.hakula@aalto.fi})} \and
Antti Rasila\thanks{Aalto University, Institute of Mathematics,
P.O. Box 11100, FI-00076 Aalto,
FINLAND ({\tt antti.rasila@iki.fi})} \and
 Matti Vuorinen\thanks{Department of Mathematics and Statistics,
FI-20014 University of Turku,
FINLAND ({\tt vuorinen@utu.fi})}}

\newcounter{minutes}
\divide\time by 60
\newcounter{hours}
\multiply\time by 60
\addtocounter{minutes}{-\time}

\begin{document}


\maketitle

\begin{abstract}
We study the problem of computing the conformal modulus of rings and quadrilaterals with strong singularities and cusps on their boundary. We reduce this problem to the numerical solution of the associated Dirichlet and Dirichlet-Neumann type boundary values problems for the Laplace equation. Several experimental results, with error estimates, are reported. In particular, we consider domains with dendrite like boundaries, in such cases where an analytic formula for the conformal modulus can be derived.  Our numerical method makes use of an $hp$-FEM algorithm, written for this very complicated geometry with strong singularities.
\end{abstract}

\begin{keywords}
conformal capacity, conformal modulus, quadrilateral modulus, $hp$-FEM, numerical conformal mapping
\end{keywords}
\begin{AMS}
65E05, 31A15, 30C85
\end{AMS}

\begin{center}
\texttt{\small FILE:~\jobname .tex 
printed: \number\year-\number\month-\number\day,
\thehours.\ifnum\theminutes<10{0}\fi\theminutes}
\end{center}

\pagestyle{myheadings}
\thispagestyle{plain}
\markboth{H. HAKULA, A. RASILA AND M. VUORINEN}{CONFORMAL MODULUS, SINGULARITIES AND CUSPS}

\section{Introduction}


The conformal modulus is an important tool in geometric function theory \cite{Ah}, and it is closely related to certain physical quantities which also occur in engineering applications. For example, the conformal modulus plays an important role in determining resistance values of integrated circuit  networks (see e.g. \cite{ps2,sl}). 
We consider both simply and doubly-connected bounded domains. By definition
such  a domain can be mapped conformally  either onto a rectangle or  onto an annulus, respectively.
For the numerical study of these two cases we define the modulus $h$ as follows.
In the simply connected case, we fix four points on the boundary  of the domain, call a domain
with these fixed boundary points a quadrilateral, and require that these four points are mapped onto vertices (0,0), (1,0), (1,h), (0,h) of the rectangle. In the doubly connected case we require
that the annulus is $\{  (x,y):  \exp(-h) < x^2+ y^2 < 1 \}\,.$
Doubly connected domains are also called ring domains or simply rings.
Surveys of the state of the art methodologies in the field are presented in the recent books by N.~Papamichael and N.~Stylianopoulos \cite{ps2} and by T.~Driscoll and L.N.~Trefethen \cite{td}. Various applications are described in \cite{hen,kuh,sl,weg}. In the past few years quadrilaterals and ring domains of increasing complexity have been studied by several authors \cite{bc,bt,dl,rcp,thg}.

We consider the problem of numerically determining the conformal modulus on certain ring 
domains with elaborate boundary. Due to the structure of the boundary, 
the problem is numerically challenging. On the other hand, the ring domain is characterized 
by a triplet of parameters $(r,m,p)$, its construction is recursive, and yet its conformal 
modulus can be explicitly given. Varying the parameter values or the recursion level of the 
construction one can increase the computational challenge and therefore this family of 
domains forms a good set of test problems. In particular, error estimates can be given.
For a figure of a domain in this family see Figure \ref{dendritefig0}. 
The  boundaries of these ring domains are point sets of dendrite type 
(i.e. continua without loops) \cite{wd}. We apply here the $hp$-FEM method developed 
in \cite{hrv1} for the computation of the moduli of these ring domains and report the
accuracy of our algorithm. 
Furthermore, we use the algorithm of \cite{hqr} to numerically approximate the canonical
conformal mapping of the above class of domains (see Figure \ref{dendritefig1} (c)).
The conjugate problem of the original ring problem solved in the approximation of the
conformal mapping
can be interpreted as a simplified crack problem.

Due to the pioneering work of I. Babuska and his coauthors \cite{bg1,bg2,bs} the optimal convergence rate
of the $hp$-FEM method is well studied and experimentally also demonstrated in some
fundamental basic situations. Our main results are the computation the moduli of
ring domains, given in the form of various error
estimates, where we compare three error norms in the case of rings with elaborate boundary:
(1) Exact error  (2) Auxiliary space estimate, based on $hp$-space theory (3) So called reciprocal
error estimate, introduced in \cite{hrv1}.
We also attain the nearly optimal convergence in accordance with the theory of \cite{bg1,bg2,bs}.
All these three error estimates behave in the same way.
Because this is the case, there is some justification to use the estimates (2) and (3) also in the
common case of applications when the exact value of the modulus is not know and hence
estimate (1) is not available. We apply our methods to compute the modulus of a quadrilateral
considered by Bergweiler and Eremenko \cite{be}. The boundary of this domain has cusp-like singularities.

\subsection{Project background and history}
This paper is a culmination of a ten-year research project, arising from questions related to the work of Betsakos, Samuelsson and Vuorinen \cite{bsv}. The original goal of the project was to develop accurate numerical tools suitable for studying effects of geometric transformations in function theory (see e.g. \cite{dv,hvv}). Experimental work towards this goal was carried out by Rasila and Vuorinen in the two small papers \cite{rv1,rv2}, and further work along the same lines was envisioned. However, it was quickly discovered that the AFEM package of Samuelsson used in the above papers is not optimal for studying very complex geometries arising from certain theoretical considerations, as the number of elements used in such computations tends to become prohibitively large. This problem led us into the higher order $hp$-FEM algorithm implemented by Hakula. Efficiency of this method for numerical computation of conformal modulus had been established in the papers \cite{hrv1} and \cite{hrv2}, the latter of which deals with unbounded domains. Recently, an implementation of this algorithms for the purpose of
numerical conformal mapping was presented
in \cite{hqr}.

\section{Preliminaries}\label{sec:preliminaries}
In this section central concepts to our discussion are introduced.
The quantities of interest from function theory are related to numerical methods, and
the error estimators arising from the basic principles are defined.
\subsection{Conformal Modulus}\label{sec:conformalmodulus}
A simply-connected domain $D$ in the complex plane $\C$ whose boundary is homeomorphic to the unit circle, is called a Jordan domain. A Jordan domain $D$, together with four distinct points $z_1, z_2, z_3, z_4$ in $\partial D\,,$ which occur in this order when traversing the boundary in the positive direction, is called a {\it quadrilateral} and denoted by $(D_1;z_1, z_2, z_3, z_4) \,.$ If $f\colon D \to fD$ is a conformal mapping onto a Jordan domain $fD$, then $f$ has a homeomorphic extension to the closure $\overline{D}$ (also denoted by $f$). We say that the {\it conformal modulus} of $(D; z_1,z_2,z_3,z_4)$ is equal to $h>0$, if there exists a conformal mapping $f$ of $\overline{D}$ onto the rectangle $[0,1]\times [0,h]$, with $f(z_1)=1+ih$, $f(z_2)=ih$, $f(z_3)=0$ and $f(z_4)=1$.

It follows immediately from the definition that the conformal modulus is invariant under conformal mappings, i.e.,
$$
\M(D; z_1,z_2,z_3,z_4) = \M(fD; f(z_1),f(z_2),f(z_3),f(z_4)),
$$
for any conformal mapping $f\colon D \to f(D)$ such that $D$ and $f(D)$ are Jordan domains.

For a curve family $\Gamma$ in the plane, we use the notation $\M(\Gamma)$ for its modulus \cite{lv}.
For instance, if $\Gamma$ is the family of all curves joining the opposite $b$-sides
within the rectangle $[0,a]\times[0,b], a,b>0,$ then $\M(\Gamma)=b/a\,.$ If we consider
the rectangle as a quadrilateral $Q$ with distinguished points $a+ib, ib,0,a$ we also have
$\M(Q;a+ib, ib,0,a)=b/a \,,$ see \cite{Ah,lv}. Given three sets $D,E,F$ we use the notation
$\Delta(E,F; D)$ for the family of all curves joining $E$ with $F$ in $D\,.$

\subsection{Modulus of a quadrilateral and Dirichlet integrals}
One can express the modulus of a quadrilateral $(D; z_1, z_2, z_3, z_4)$
in terms of the solution of the Dirichlet-Neumann problem as follows.
Let $\gamma_j$, $j=1,2,3,4$ be the arcs of
$\partial D$ between $(z_4, z_1)\,,$ $(z_1, z_2)\,,$ $(z_2, z_3)\,,$
$(z_3, z_4),$ respectively. If $u$ is the (unique) harmonic solution of
the Dirichlet-Neumann problem with boundary values of $u$ equal to $0$ on
$\gamma_2$, equal to $1$ on $\gamma_4$ and with $\partial u/\partial n =
0$ on $\gamma_1 \cup \gamma_3\,,$ then by \cite[p. 65/Thm 4.5]{Ah}:
\begin{equation}
\label{qmod}
\M(D;z_1,z_2,z_3,z_4)=
\iint_D |\nabla u|^2\,dx\,dy.
\end{equation}
The function $u$ satisfying the above boundary conditions is called the
{\it potential function} of the quadrilateral $(D;z_1, z_2, z_3, z_4)$.

\subsection{Modulus of a ring domain and Dirichlet integrals}
Let $E$ and $F$ be two disjoint compact sets in the extended
complex plane ${\C_\infty}$. Then one of the sets $E,$ $F$ is bounded and
without loss of generality we may assume that it is $E\,.$ If both $E$ and $F$ are connected
and the set $R={\C_\infty} \setminus(E \cup F)$ is connected,
then $R$ is called a {\it ring domain}. In this case $R$ is a doubly
connected plane domain. The {\it capacity} of $R$ is defined by
$$
\Cap R=\inf_u \iint_D |\nabla u|^2\,dx\,dy,
$$
where the infimum is taken over all nonnegative, piecewise
differentiable functions $u$
with compact support in $R\cup E$ such that $u=1$ on $E$.
It is well-known that there exists a unique harmonic function on $R$ with
boundary values $1$ on $E$ and $0$ on $F$. This function is called the potential
function of the ring domain $R$, and it minimizes the above integral. In other words,
the minimizer may be found by solving the Dirichlet problem for the Laplace equation in $R$ with
boundary values $1$ on the bounded boundary component $E$ and $0$ on the
other boundary component $ F\,.$
A ring domain $R$ can be
mapped conformally onto the annulus $\{z:e^{-M}<|z|<1\}$, where
$M=\M(R)$ is the conformal modulus
of the ring domain $R\,.$ The modulus and capacity of a ring
domain are connected by the simple identity
$\M(R)=2\pi/\Cap R$.
For more information on the modulus of a ring domain
and its applications in complex analysis the reader is
referred to \cite{Ah,hen,kuh,ps2}.

\subsection{Hyperbolic Metrics}
The hyperbolic geometry in the unit disk is a powerful tool of classical complex analysis.
 We shall now briefly review some of the main features of this geometry, necessary for what follows.
First of all, the hyperbolic distance
between  $x,y \in {\mathbb{D} }$ is given by
$$
\rho_{ \mathbb{D} } (x,y) = 2\, {\rm arsinh \left(\frac{|x-y|}{
\sqrt{(1-|x|^2)(1- |y|^2)}} \right)}.
$$
In addition to the unit disk ${\mathbb{D} }$, one  usually also studies the upper half plane $\IH$ as a model of the hyperbolic  geometry.
For $x,y\in \IH$ we have  ($x= (x_1, x_2)$)
$$
\rho_{  \IH}(x,y) = {\rm arcosh} \bigg(1 + \frac{|x-y|^2}{2 x_2 
 y_2}\bigg) \,.
$$
If there is no danger of confusion,
we denote both $\rho_{\IH}(z,w)$ and $\rho_{\ID}(z,w)$
simply by $\rho(z,w)$. 
We assume that the reader is familiar with some basic facts about these geometries: geodesics, hyperbolic length minimizing curves, are circular arcs orthogonal to the boundary in each case.

Let $z_1,z_2,z_3,z_4$ be distinct points in $\mathbb{C}$. We define the {\it absolute
(cross) ratio} by
\begin{equation}
\label{cross}
|z_1,z_2,z_3,z_4| = \frac{|z_1-z_3|\,|z_2-z_4|}{|z_1-z_2|\,|z_3-z_4|}.
\end{equation}
This definition can be extended for $z_1,z_2,z_3,z_4\in \mathbb{C}_\infty$ by taking the limit.
An important property of M\"obius transformations is that they preserve
the absolute ratios, i.e.
$$
|f(z_1),f(z_2),f(z_3),f(z_4)|= |z_1,z_2,z_3,z_4|,
$$
if $f\colon \mathbb{C}_\infty\to \mathbb{C}_\infty$ is a M\"obius transformation. In fact, a mapping
$f\colon \mathbb{C}_\infty\to \mathbb{C}_\infty$ is a M\"obius transformation if and only if $f$ it is sense-preserving and
preserves all absolute ratios.

Both for $(\ID,\rho_{\ID})$ and $(\IH,\rho_{\IH})$ one can define the hyperbolic distance in terms of the absolute ratio. Since the absolute ratio is invariant under M\"obius transformations, the hyperbolic metric also remains invariant under these transformations. In particular, any M\"obius transformation of   ${\mathbb{D} }$ onto $\IH$ preserves the hyperbolic distances. A standard reference on hyperbolic metrics is \cite{Beardon}.

\subsection{$hp$-FEM}\label{sec:hpfem}

In this work the natural quantity of interest is always related to the Dirichlet energy.
Of course, the finite
element method (FEM) is an energy minimizing method and therefore an obvious choice.
The continuous Galerkin $hp$-FEM algorithm used throughout this paper is 
based on our earlier work 
\cite{hrv1}. Brief outline of the relevant features used in numerical examples below is:
Babu\v{s}ka-Szabo -type $p$-elements,
curved elements with blending-function mapping for exact geometry,
rule-based meshing for geometrically graded meshes, and
in the case of isotropic $p$ distribution, hierarchical solution for all $p$.
The main new feature considered here is the introduction of auxiliary subspace techniques
for error estimation.

For the types of problems considered here, theoretically optimal conforming $hp$-adaptivity 
is hard. The main
difficulty lies in mesh adaptation since the desired geometric or exponential grading is
not supported by standard data structures such as Delaunay triangulations.
Thus, the approach advocated here is a hybrid one, where the problem is first solved
using an \textit{a priori} $hp$-algorithm after which the quality of the solution is
estimated using error estimators specific both for the problem and the method, provided the
latter are available. For instance, the exact solution or
for problems concerning the conformal modulus
the so-called reciprocal error estimator. The a priori algorithm is modified if
the error indicators suggest modifications. If this occurs, the solution process is started
anew.

In the numerical examples below the computed results are
measured with both kinds of error estimators giving us high confidence in the validity 
of the results and the chosen methodology.

\subsubsection{Auxiliary Subspace Techniques}\label{sec:errorestimator}
Consider the 
abstract problem setting with $V$ as the standard piecewise polynomial
finite element space on some discretization $T$ of the computational 
domain $D$. Assuming that the exact solution $u \in H_0^1(D)$ has finite energy,
we arrive at the approximation problem:
Find $\hat{u} \in V$ such that
\begin{equation}\label{eq:approximation}
  a(\hat{u},v)=l(v)\ (= a(u,v)),\ \forall v \in V,
\end{equation}
where $a(\cdot,\cdot)$ and $l(\cdot)$, are the bilinear form and the load
potential, respectively.
Additional degrees of freedom can be introduced by enriching the space $V$.
This is accomplished via introduction of an auxiliary subspace or ``error space'' 
$W \subset H_0^1(D)$
such that $V \cap W = \{0\}$.
We can then define the error problem: Find $\epsilon \in V$ such that
\begin{equation}\label{eq:error}
  a(\epsilon,v)=l(v)- a(\hat{u},v) (= a(u-\hat{u},v)),\ \forall v \in W.
\end{equation}
In 2D the space $W$, that is, the additional unknowns, can be associated 
with element edges and interiors. Thus, for $hp$-methods this kind of
error estimation is natural. For more details on optimal selection of
auxiliary spaces, see \cite{hno}.

The solution $\epsilon$ of (\ref{eq:error}) is called the \textit{error function}. 
It has many useful properties for both theoretical and practical considerations.
In particular, the error function can be numerically evaluated and
analysed for any finite element solution. This property will be used in the following.
By construction, the error function is identically zero at the mesh points.
In Figure~\ref{dendritefig2} one instance of a contour plot of the error function 
(with a detail) is shown. This gives an excellent way to get a qualitative view
of the solution which can be used to refine the discretization in the $hp$-sense.

Let us denote the error indicator by a pair $(e,b)$, where $e$ and $b$ refer to
added polynomial degrees on edges and element interiors, respectively.
It is important to notice that the estimator requires a solution of a linear system.
Assuming that the enrichment is fixed over the set of $p$ problems, it is clear that
the error indicator is expensive for small values of $p$ but becomes asymptotically
less expensive as the value of $p$ increases.
Following the recommendation of \cite{hno}, our choice in the sequel is
$(e,b)=(1,2)$ unless specified otherwise.

\begin{remark}
In the case of $(0,b)$-type or pure bubble indicators, the system is not connected and
the elemental error indicators can be computed independently, and thus in parallel.
Therefore in practical cases one is always interested in relative performance
of $(0,b)$-type indicators.
\end{remark}

\subsection{Reciprocal Identity and Error Estimation}
Let $Q$ be a quadrilateral defined by points $z_1,z_2,z_3,z_4$ and boundary curves as
in Section \ref{sec:conformalmodulus} above.
The following reciprocal identity holds:
\begin{equation} \label{recipridty}
\M(Q; z_1,z_2,z_3,z_4) \M(Q; z_2,z_3,z_4, z_1) =1.
\end{equation}
As in \cite{hrv1,hrv2}, we shall use the test functional
\begin{equation}
\label{recitest}
\big| \M(Q; z_1,z_2,z_3,z_4) \M(Q; z_2,z_3,z_4, z_1) -1\big|
\end{equation}
which by \eqref{recipridty} vanishes identically, as an error estimate.

As noted above, the error function $\epsilon$ can be analysed in the sense of FEM-solutions.
Our goal is to relate the error function given by auxiliary space techniques and the 
reciprocal identity arising naturally from the geometry of the problem.
Let us first define the energy of the error function $\epsilon$ as
\begin{equation}\label{eq:errorenergy}
  \mathcal{E}(\epsilon) = \iint_D |\nabla \epsilon|^2\,dx\,dy.
\end{equation}

Using (\ref{eq:errorenergy}) the reciprocal error estimation and the error function $\epsilon$ introduced above 
can be connected as follows:
Let $a_1$ and $a_2$ be the moduli of the original and conjugate problems, $\epsilon_1$, $\epsilon_2$, and $\hat{\epsilon}_1 =  \mathcal{E}(\epsilon_1)$, $\hat{\epsilon}_2 =  \mathcal{E}(\epsilon_2)$, the errors and their energies, respectively.
Taking $\hat{\epsilon} =  \max\{|\hat{\epsilon}_1|,|\hat{\epsilon}_2|\}$ we get via direct computation:
\begin{equation} 
|1-(a_1 + \hat{\epsilon}_1)(a_2 + \hat{\epsilon}_2)| \leq 
  |a_1\hat{\epsilon}_2 + a_2\hat{\epsilon}_1 + \hat{\epsilon}_1 \hat{\epsilon}_2|
  \leq 2\, \hat{\epsilon} \max\{a_1,1/a_1\}+O(\hat{\epsilon}^2).
\end{equation}
Neglecting the higher order term one can solve for $\hat{\epsilon}$ and compare this
with the estimates given by the individual error functions.
\section{The Dendrite}

A compact connected set in the plane is called dendrite-like if it contains no loops. We introduce a new parametrized family of ring domains whose boundaries have dendrite-like boundary structure and whose modulus is explicitly known in terms of parameters.  In numerical conformal mapping one usually considers domains whose
boundaries consist of finitely many piecewise smooth curves. Very
recently in \cite{rcp} these authors considered conformal mapping
onto domains whose boundaries have ``infinitely many sides'', i.e.,
are obtained as result of a recursive construction. One of the
examples considered in \cite{rcp} was the domain whose boundary
was the von Koch snowflake curve.

In this section we will give a construction of a ring domain
whose complementary components are ${\mathbb C} \setminus \mathbb{D}$
and a compact connected subset $C(r,p,m)$ of the unit disk $\mathbb{D}$
depending on two positive integer parameters $p,m$ and a real number $r>0$.
The set $C(r,p,m)$ consists of finitely many pieces, each of which is a smooth curve, and the set
is acyclic, i.e., does not contain any loops. The number of
pieces is controlled by the integers $(p,m)\,$ and can be arbitrarily
large when $p$ and $m$ increase.

\subsection{Theory}


Recall that the Gr\"otzsch ring $R_G(r)=\mathbb{D} \setminus [0,r]$,
$r\in(0,1)$ has the capacity $\Cap(R_G(r))=2\pi/\mu(r)$, were $\mu(r)$
is the Gr\"otzsch modulus function (cf. \cite[Chapter 5]{avv}):
$$
\mu(r) = \frac{\pi}{2}\frac{\K(r')}{\K(r)} \, , \text{ and }
\K(r)=\int_0^1 \frac{dx}{\sqrt{(1-x^2)(1-r^2 x^2)}}\, ,
$$
with usual notation $r'=\sqrt{1-r^2}$.
Let $r\in(0,1)$, and let $D_r= \mathbb{D} \setminus
\big([-r,r]\cup[-ir,ir]\big)$. The conformal mapping 
$f(z)=\sqrt[4]{-z}$ maps the Gr\"otzsch ring $R_G(r)$ (excluding the positive) real axis onto the sector $\{z: |\arg z|<\pi/2\}$. Let $u_r$ be the potential function associated with $R_G(r)$. Then, by symmetry, it follows that the potential function $u_r\circ f$ can be extended to the domain $D_r$ by Schwarz symmetries so that it solves the Dirichlet problem associated with the conformal capacity of $D_r$. It follows that $\Cap(D_r) = 8\pi/\mu(r^4)$. Obviously, the similar construction is possible for any integer $m\ge 3$.

One may continue the process to obtain further generalizations. Start with a generalized Gr\"otzsch
ring with $m\geq 3$ branches. Choose one of the vertices of the
interior component. Map this point to the origin by a M\"obius
automorphism of the unit disk. Make a branch of degree $p\geq 2$ to the
origin by using the mapping $z\mapsto z^{1/p}$, and extend the potential function to the whole disk by using Schwarz symmetries. The resulting ring has capacity $2\pi mp/\mu(r^m)$. An example of the construction is given in
Figure~\ref{dendritefig0}.

Again, it is possible to further iterate the above construction to obtain ring domains with arbitrarily complex dendrite-like boundaries. Let $m\ge 3$, $M\ge 1$, and let be integers such that $p_j\ge 2$ for all $j=1,2,\ldots,M$. For each $j=1,2,\ldots,M$, choose one of the vertices $z_j$ of the interior component. Let $w_j$ be the point on the line $tz_j$, $t>0$ so that $|w_j|=1$. We may assume that $z_j<0$, $w_j=-1$ and the line segment $[-1,z_j]$ does not intersect with the interior component except at the point $z_j$. Map the point $z_j$ to the origin by M\"obius automorphism $g_j$ of the unit disk so that $g_j(-1)=-1$. Now map the domain $\mathbb{D}\setminus [-1,0]$ onto the symmetric disk sector by the mapping $h_j(z)=z^{1/p_j}$, and extend the potential function to the whole disk by using Schwarz symmetries. By repeating this construction for all $j=1,2,\ldots,M$, we obtain a ring domain with the conformal capacity 
\[
\frac{2\pi mp_1p_2\cdots p_M}{\mu(r^m)}.
\]

\subsection{Numerical Experiments}\label{sec:dendriteexperiments}

\begin{figure}
\centering
\subfloat[Generalized Gr\"otzsch ring $R_G(r, m)$: $r~=~1/4$, $m=6$.]
{\includegraphics[width=0.45\textwidth]{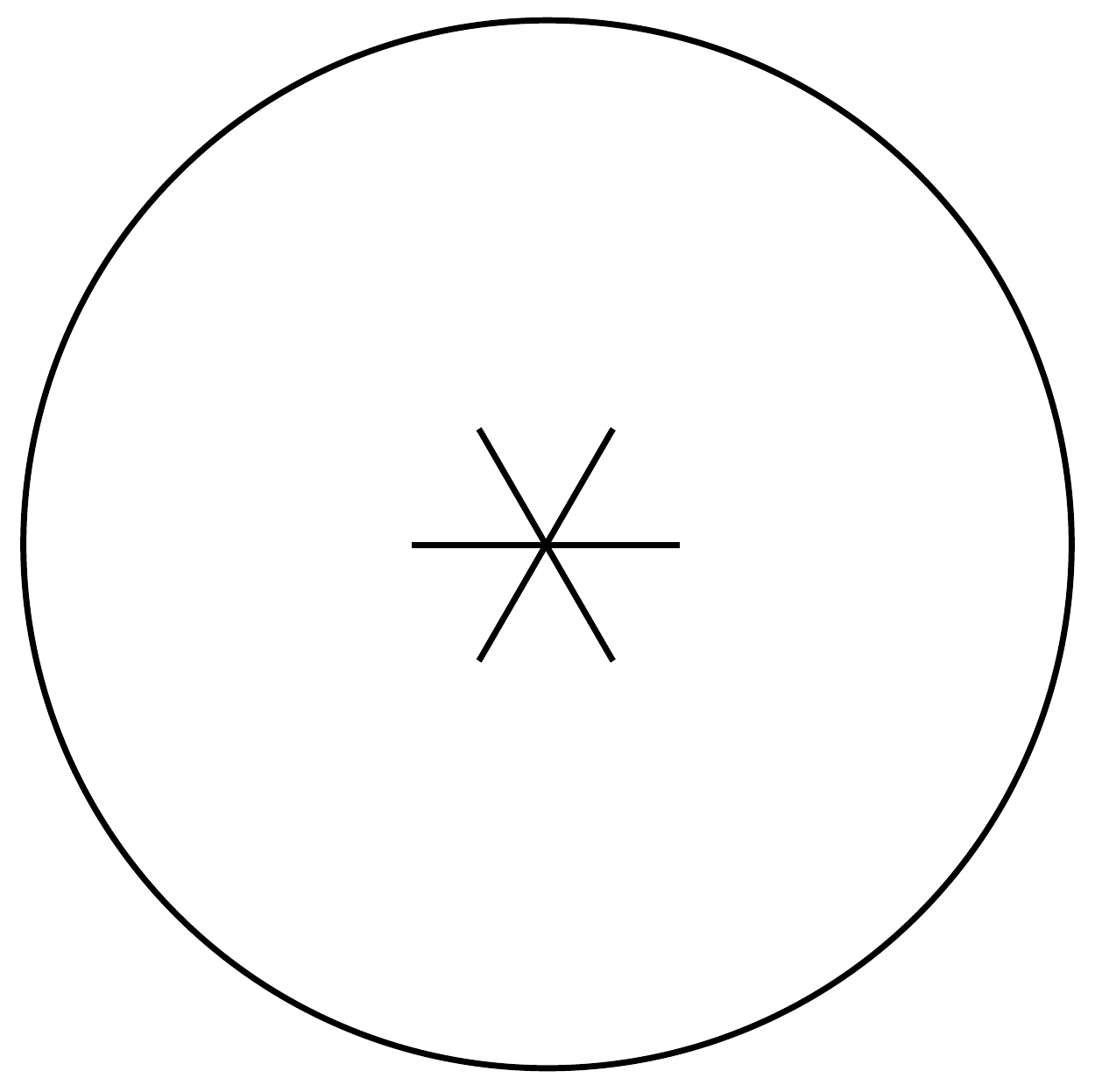}}\quad
\subfloat[Map the chosen point to the origin by a M\"obius
automorphism of the unit disk.]
{\includegraphics[width=0.45\textwidth]{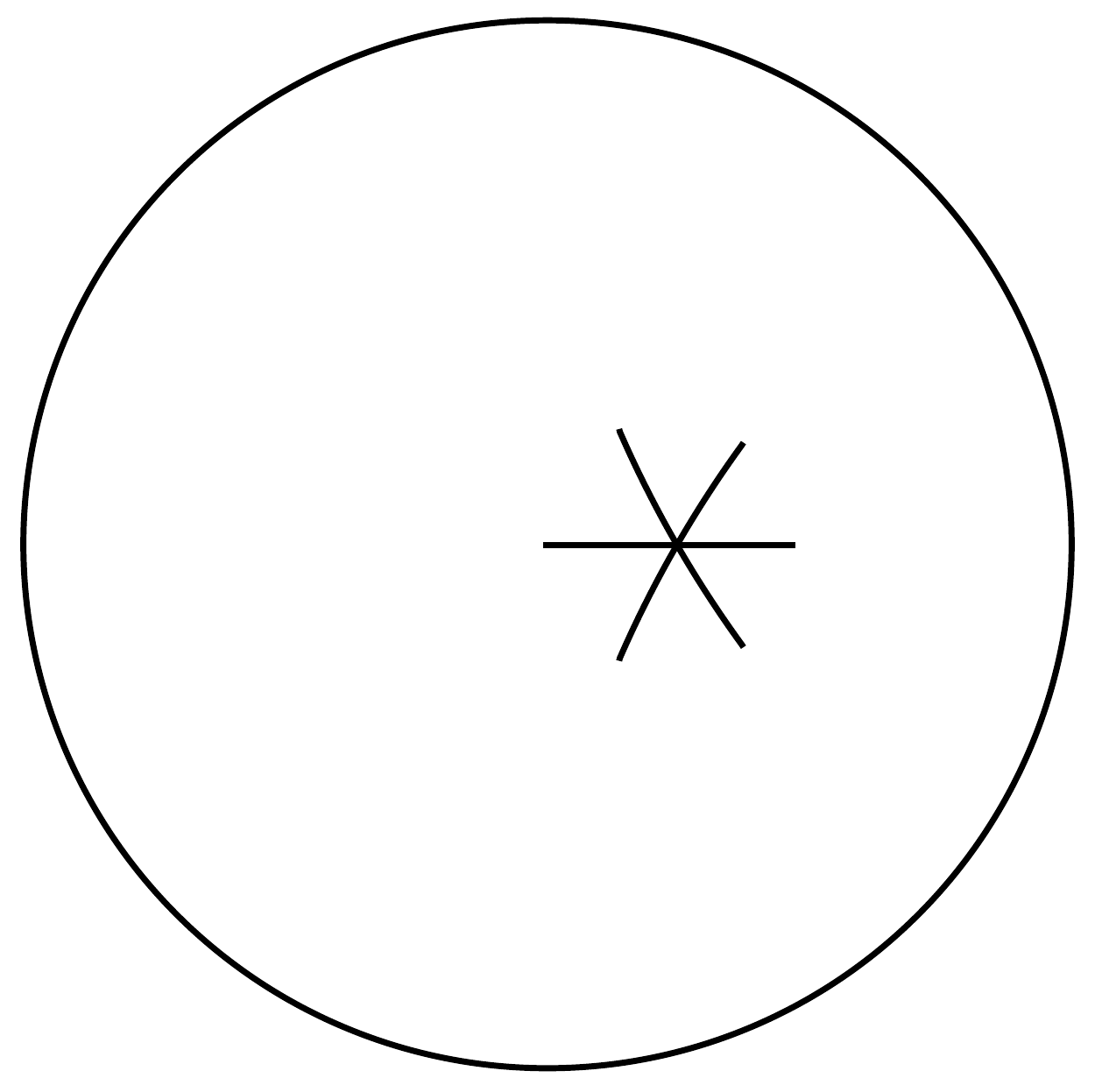}}\\
\subfloat[Make a branch of degree $p=7$.]
{\includegraphics[width=0.45\textwidth]{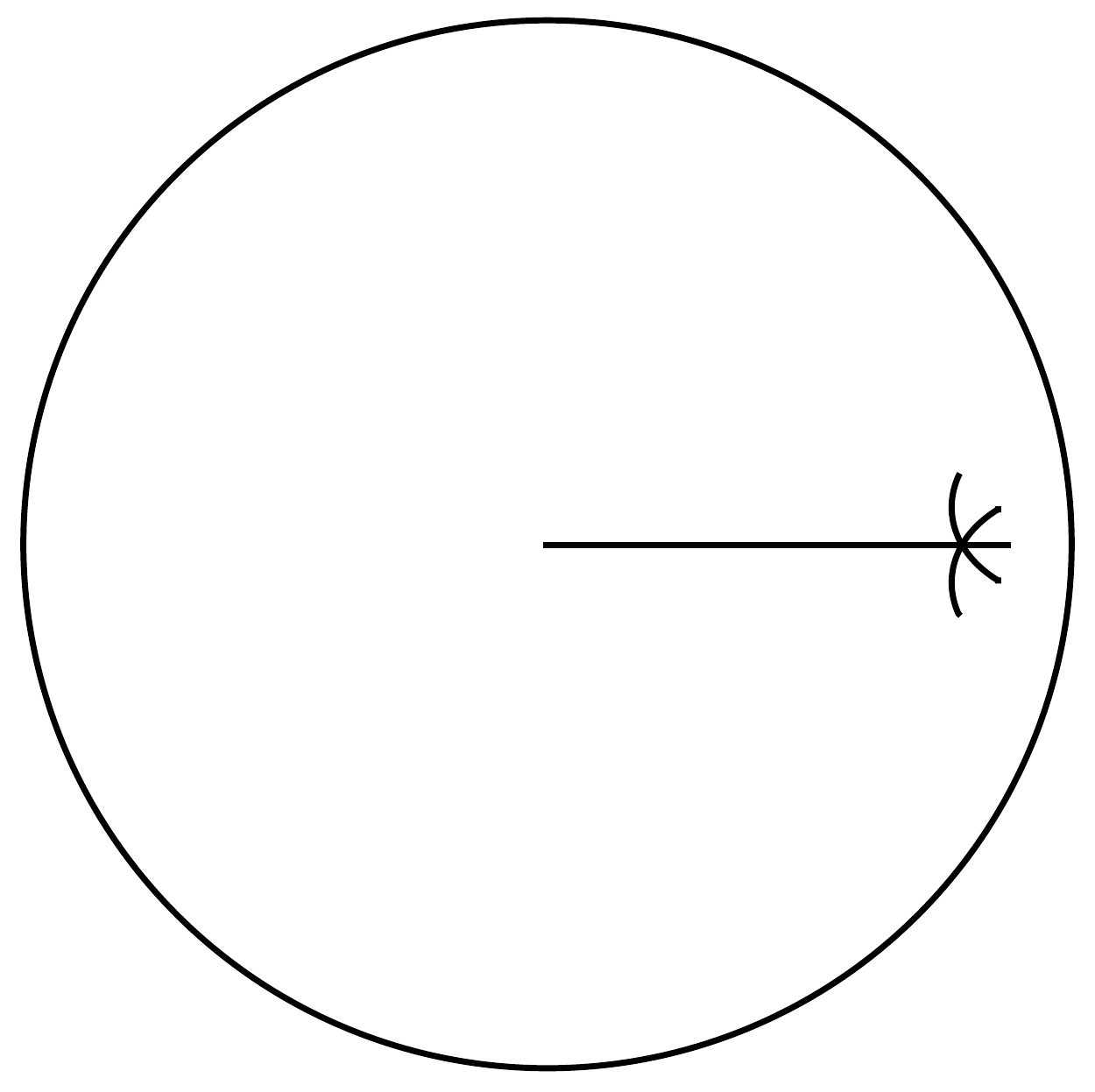}}\quad
\subfloat[Extend the potential function to the whole disk.]
{\includegraphics[width=0.45\textwidth]{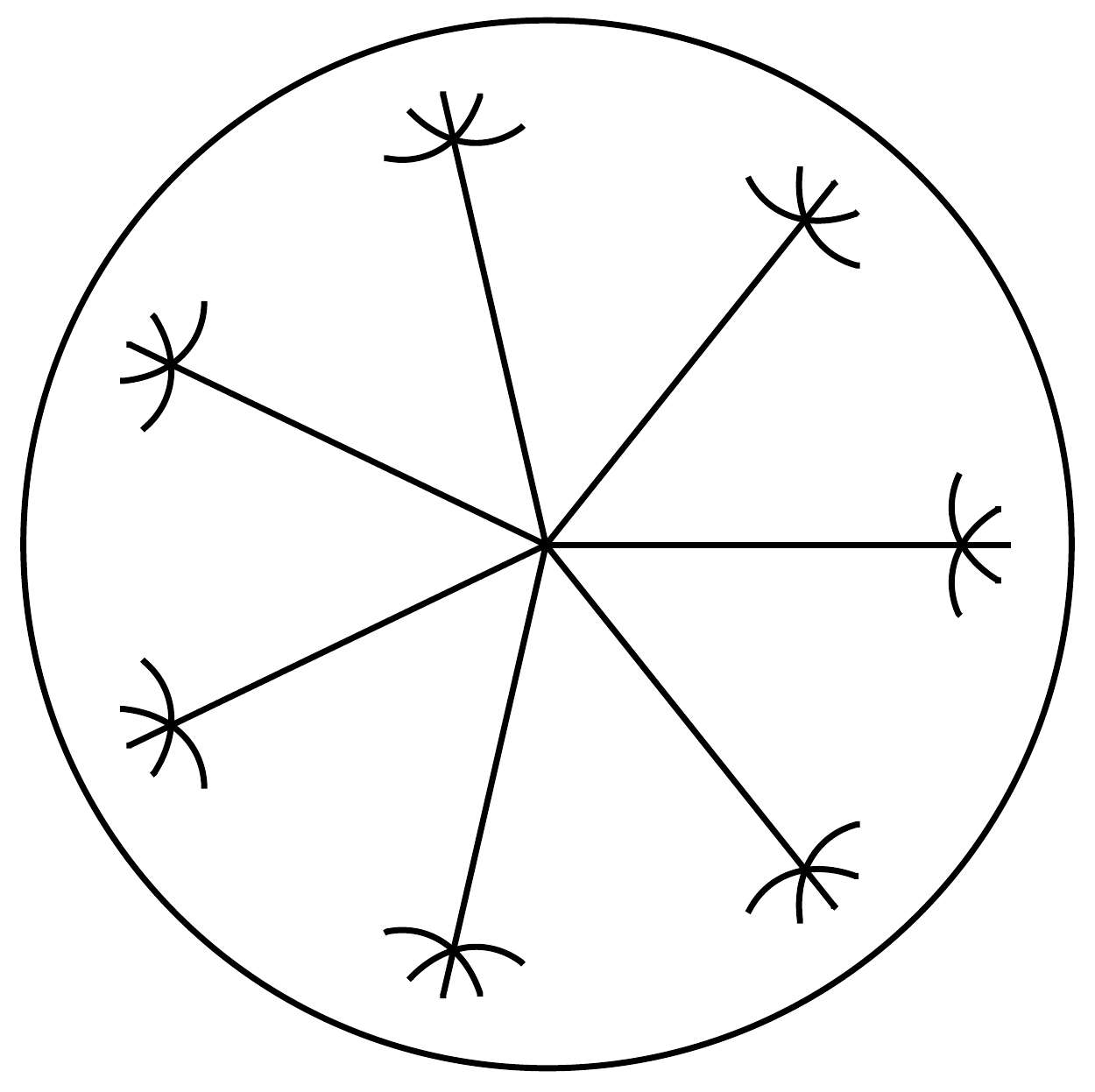}}
\caption{Dendrite Construction: $r=1/4$, $m=6$, $p=7$.}\label{dendritefig0}
\end{figure}

\begin{figure}
\centering
\subfloat[Domain.]{\includegraphics[width=0.30\textwidth]{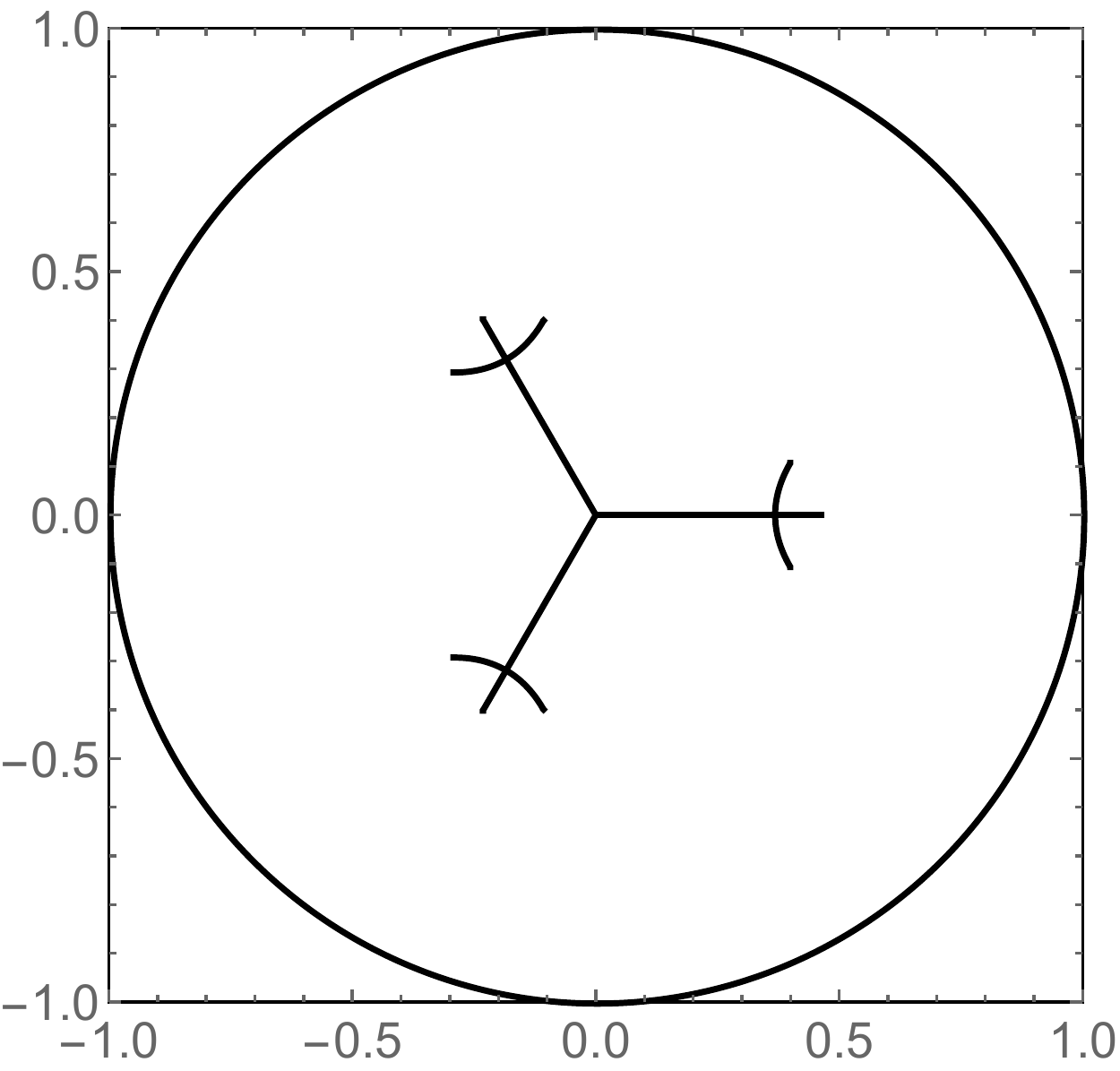}}\quad
\subfloat[Potential.]{\includegraphics[width=0.30\textwidth]{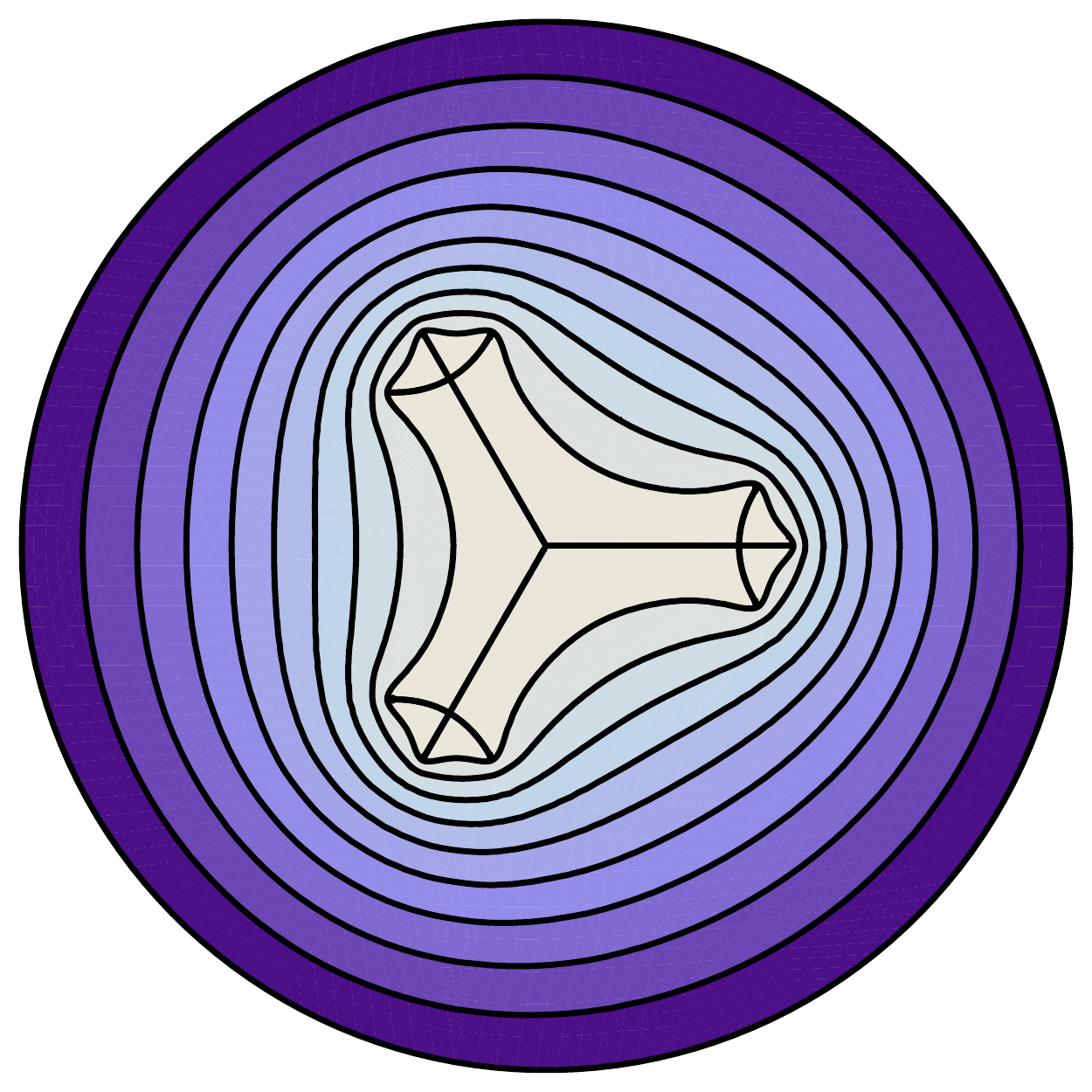}}\quad
\subfloat[Conformal mapping.]{\includegraphics[width=0.30\textwidth]{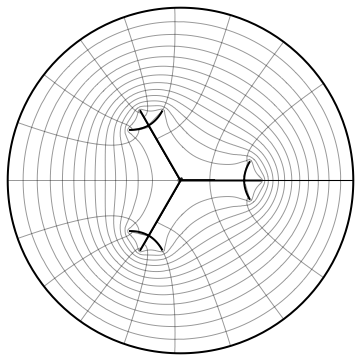}}
\caption{Dendrite 1: $r=1/20$, $m=4$, $p=3$.}\label{dendritefig1}
\end{figure}

\begin{figure}
\centering
\subfloat[Domain.]{\includegraphics[width=0.30\textwidth]{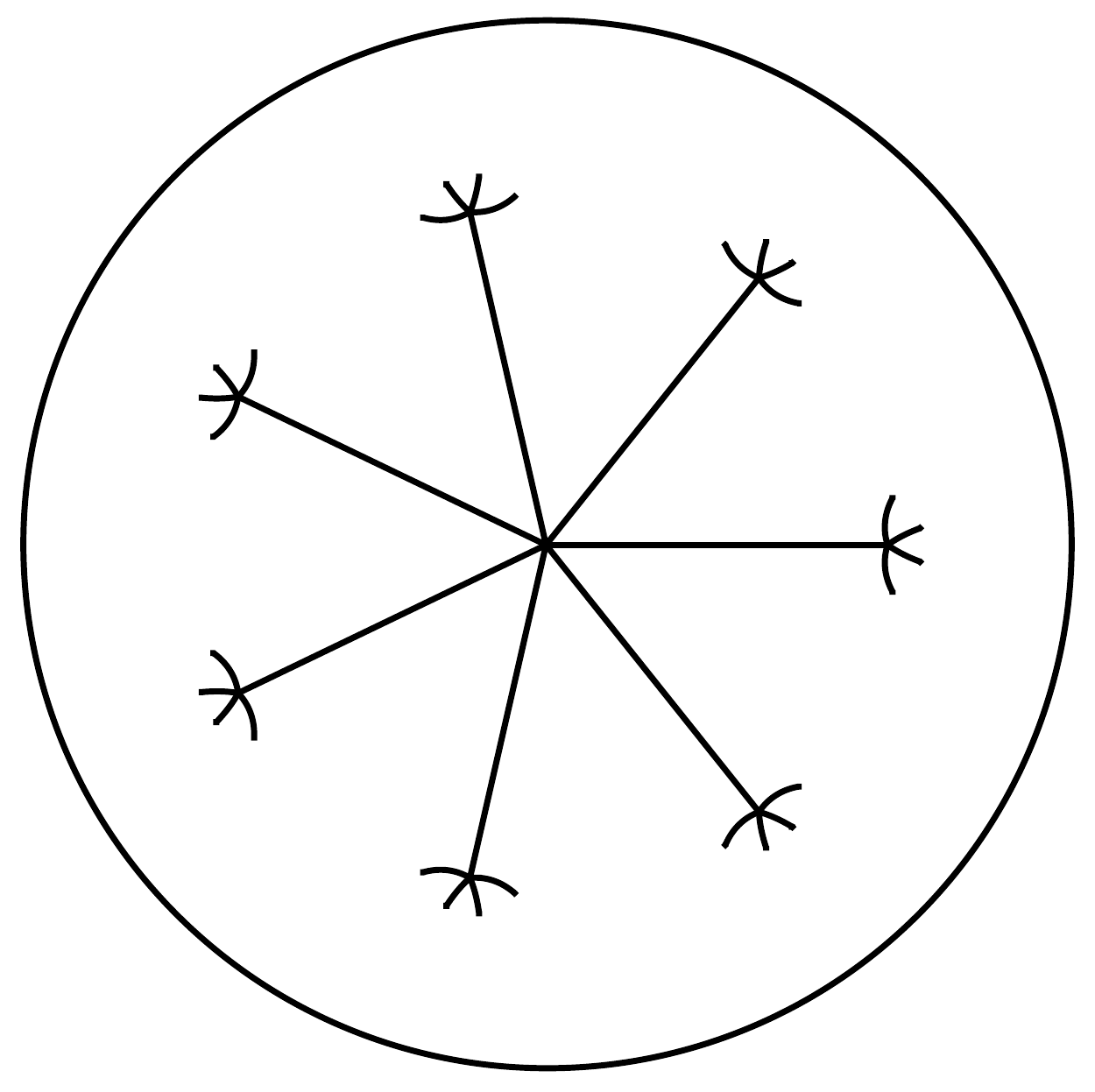}}\quad
\subfloat[Contour plot: Error function $(e,b)=(1,2)$, $hp,p=10$.]{\includegraphics[width=0.30\textwidth]{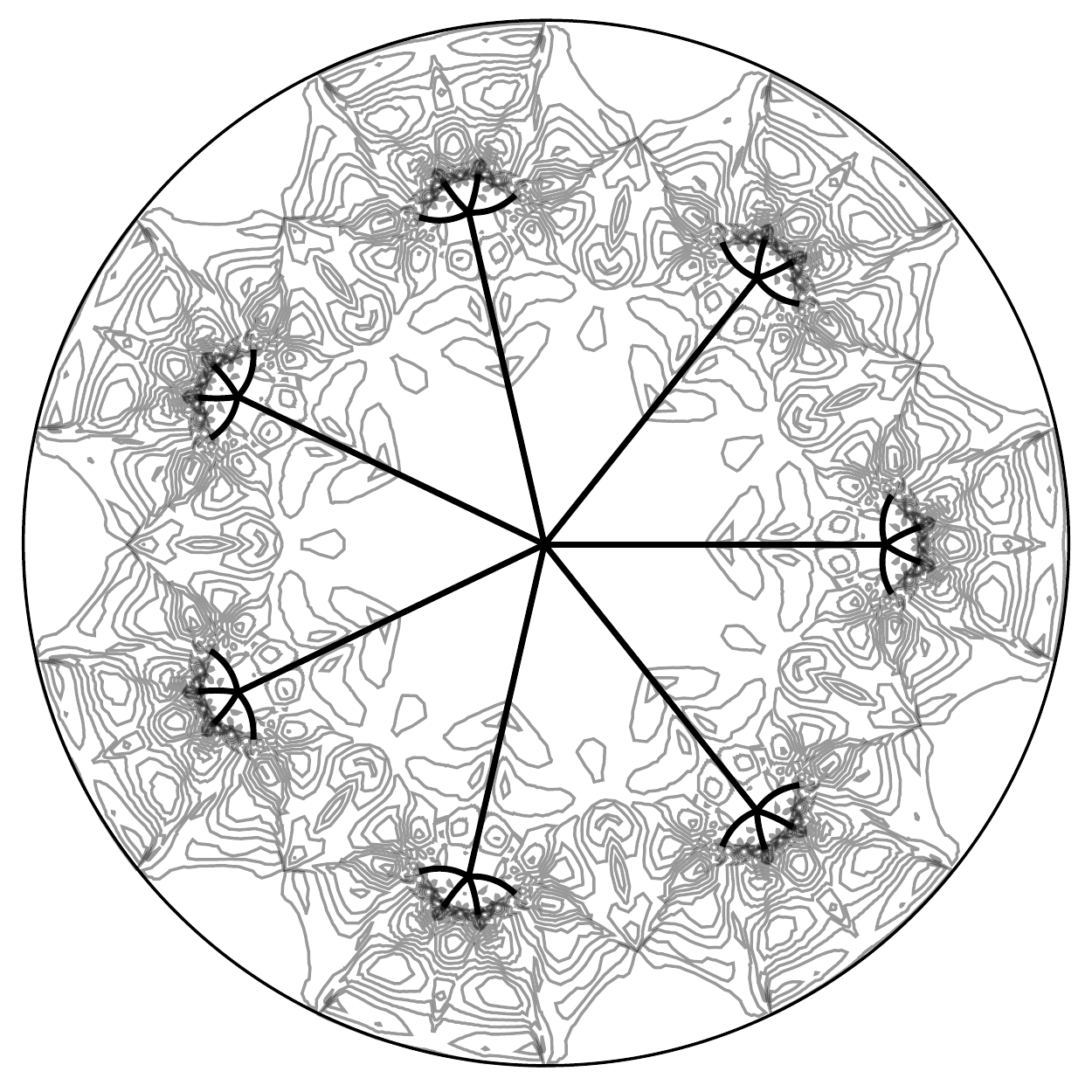}}\quad
\subfloat[[Contour plot: Error function detail.]%
{\includegraphics[width=0.30\textwidth]{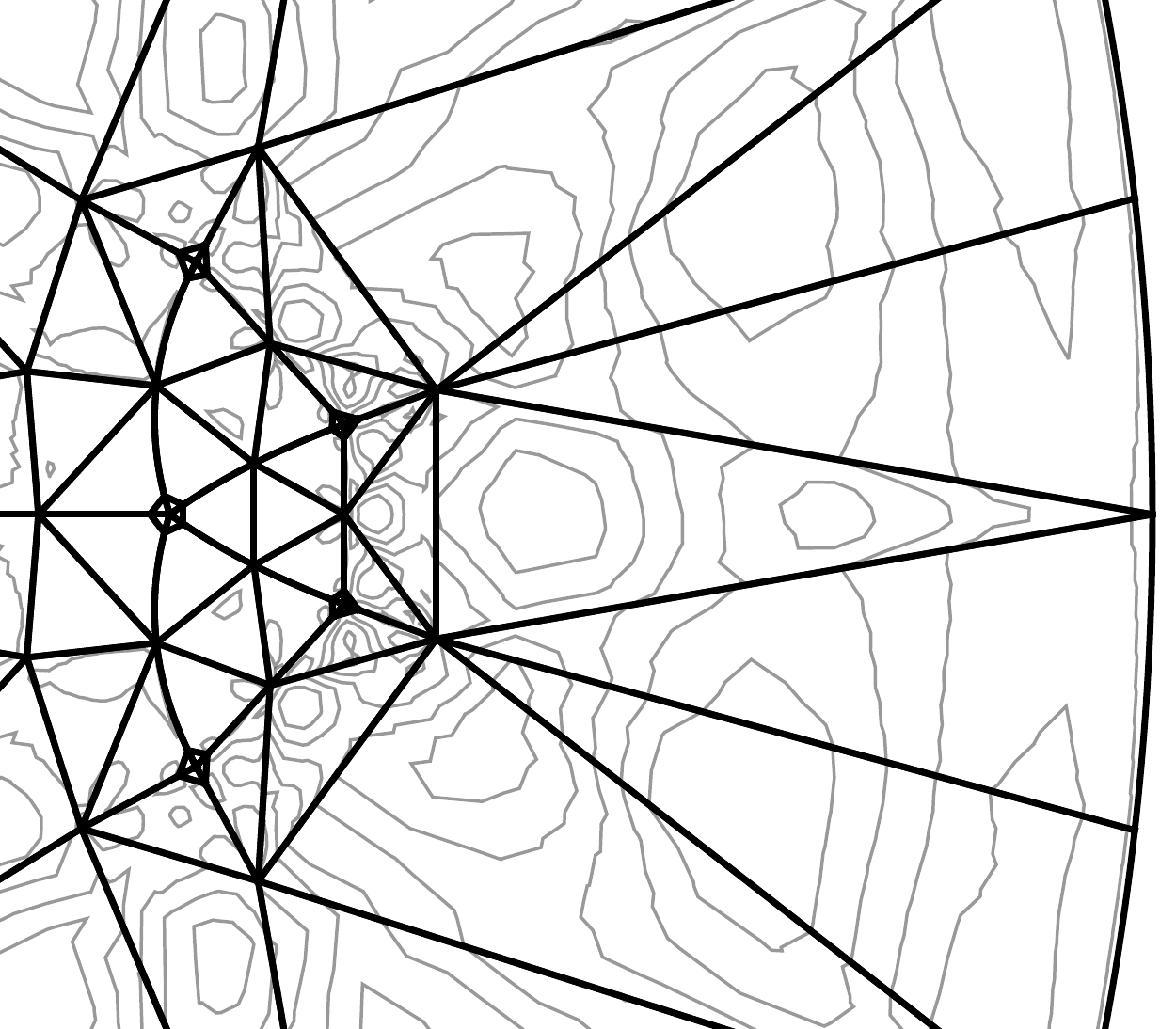}}
\caption{Dendrite 2: $r=1/20$, $m=5$, $p=7$.}\label{dendritefig2}
\end{figure}

\begin{figure}
\centering
\includegraphics[width=0.45\textwidth]{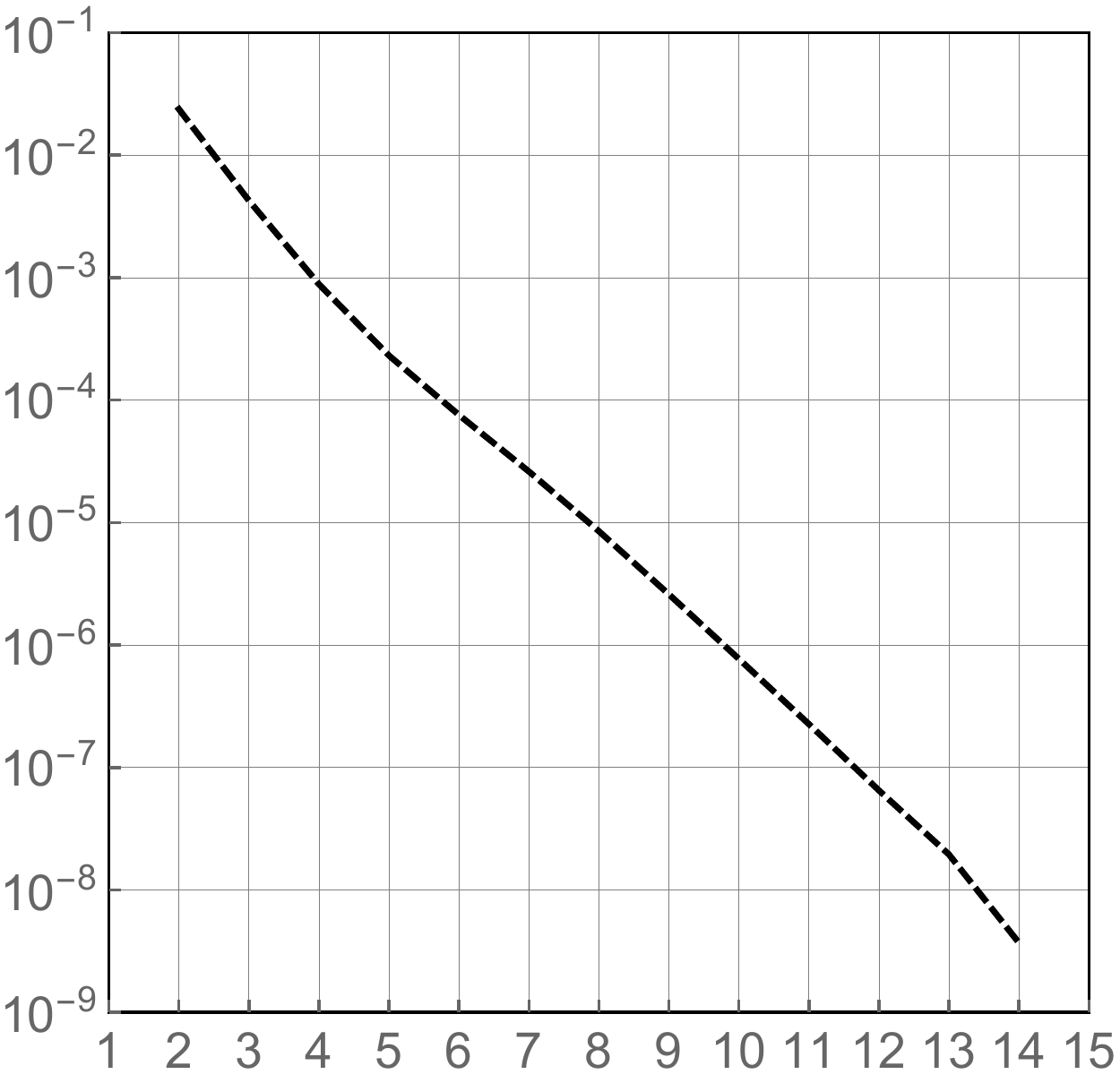}
\caption{Dendrite 1: Reciprocal error; $\log$-plot: Error vs $p$.}
\label{fig:d1recip}
\end{figure}

We consider two cases described in Table~\ref{dendritefig0} and Figures~\ref{dendritefig1}~and~\ref{dendritefig2}. 
Using the $hp$-refinement
strategy at the tips of the dendrite and the inner angles ($120^{\circ}$) we obtain exponential
convergence in the reciprocal error (Figure~\ref{fig:d1recip}).
In Figure~\ref{dendritefig2} we also show the error function of type (1,2) over the whole domain
as well as a detail which clearly shows the non-locality of the error function.
As expected the errors are concentrated at the singularities and in the elements
connecting the singularities to the boundary. One should bear in mind, however, that the
reciprocal errors are small already at $p=10$ used in the figures.

\begin{figure}
\centering
\subfloat[Estimated error.]{\includegraphics[width=0.45\textwidth]{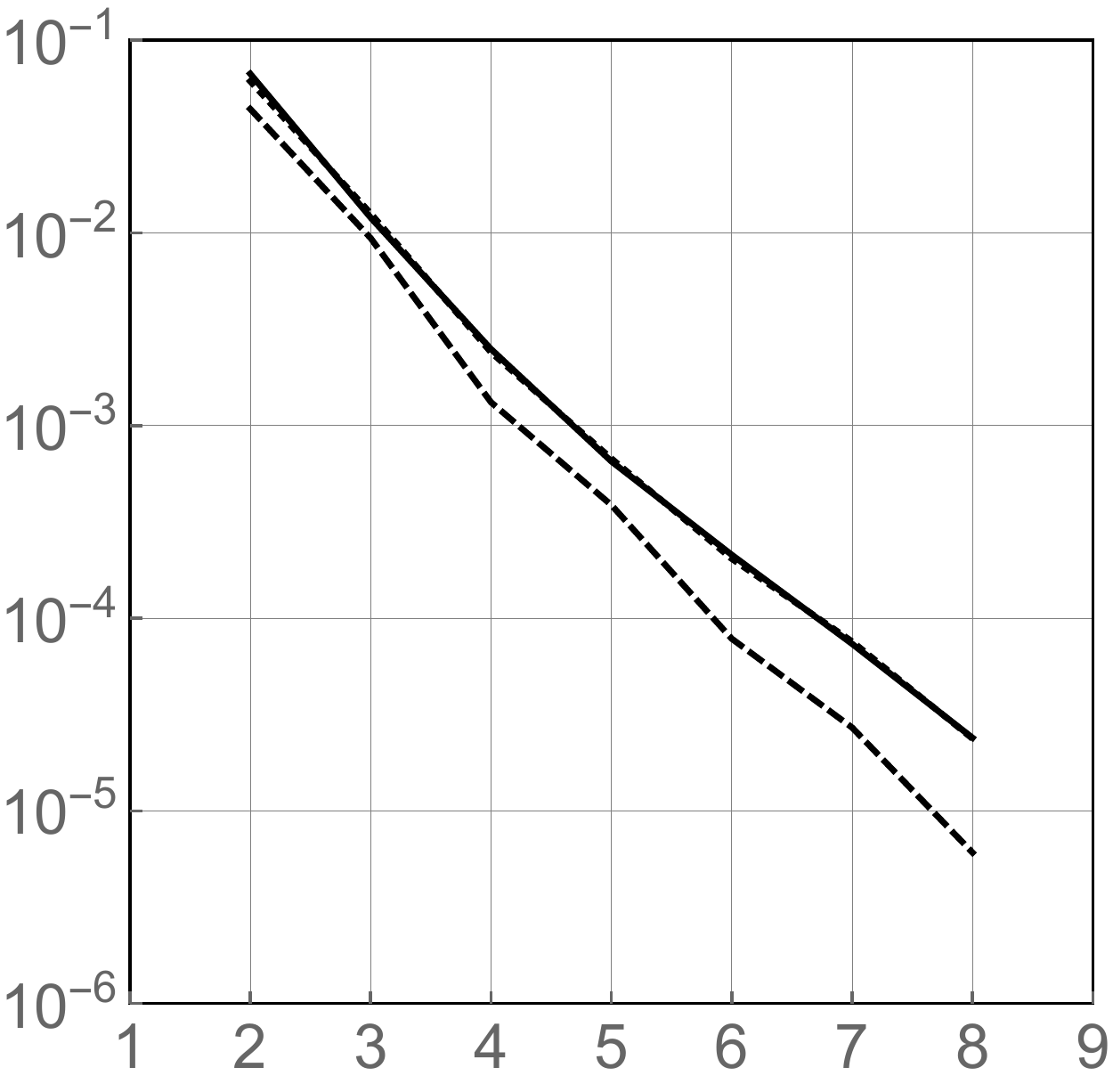}}\quad
\subfloat[Estimated error (Conjugate).]{\includegraphics[width=0.45\textwidth]{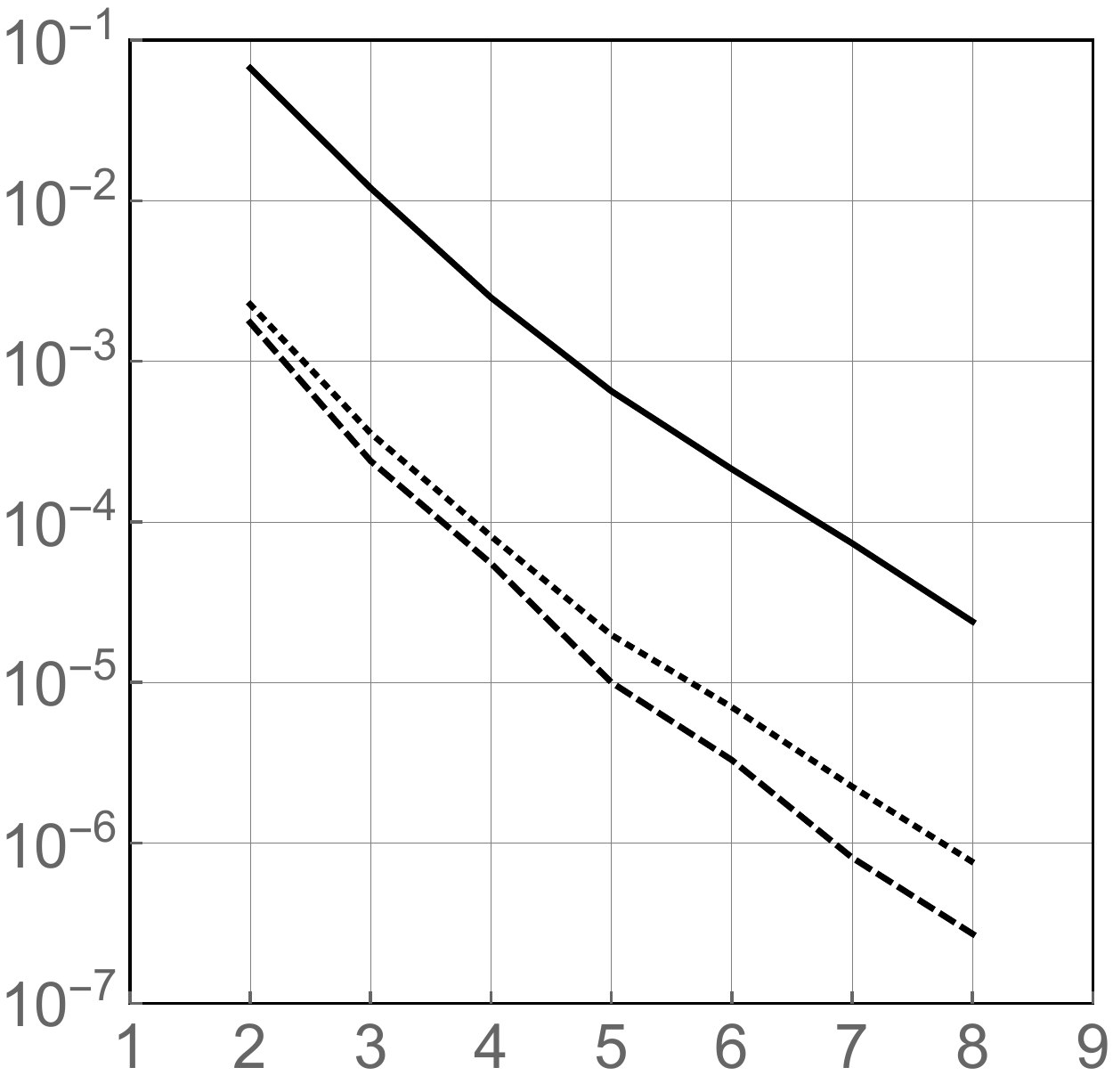}}
\caption{Dendrite 1: Estimated errors; $\log$-plot: Error vs $p$; Solid line = Reciprocal estimate, Dashed line = Auxiliary space estimate, Dotted line = Exact error.}\label{dendritefig4}
\end{figure}

The error estimates are shown in Figure~\ref{dendritefig4}.
The effect of error balancing is evident. For the larger capacity, the reciprocal error
conincides with the true error, but overestimates the smaller one. However, in both cases
the rates are correct, only the constant is overly pessimistic. The auxiliary space
error estimate underestimates the error slightly, again with the correct rate. This is exactly
what we would expect, since increasing the polynomial orders in the scheme should increase the
error estimate, if the underlying solution has converged to the correct solution.

\begin{table}[ht]
\caption{
  Tests on dendrite-problems. The errors are given as {$|\lceil\log_{10}|\mathrm{error}|\rceil|$}.
}\label{tbl:dendrite}
\renewcommand\arraystretch{1}
\noindent

\begin{center}\footnotesize
\begin{tabular}{|l|l|l|l|l|l|l|}
\hline
Case & Parameters & Method & Errors & Sizes & $M(Q_s)$\\
\hline
1 & $r=1/20$, $m=4$, $p=3$ & $hp$, $p=16$ & 9\ (9) & 159865\ (160161) & 5.63968609980242\\
2 & $r=1/20$, $m=5$, $p=7$ & $hp$, $p=12$ & 9\ (9) & 199921\ (200809) & 13.437951766839522\\
\hline
\end{tabular}
\end{center}
\end{table}

\section{Domains with Cusps}


\subsection{Hyperbolic Quadrilateral}
\label{sec:hyperbolicquadrilateral}

Let $Q_s$ be the quadrilateral whose sides are circular arcs
perpendicular to the unit circle with vertices $e^{is}$, $e^{(\pi-s)i}$,
$e^{(s-\pi)i}$ and $e^{-si}$. We call quadrilaterals of this type hypebolic quadrilaterals as their sides are geodesics in the hyperbolic geometry of the unit disk. We approximate values of the modulus of $Q_s$.

Next we determine a lower bound for the modulus of a hyperbolic quadrilateral. 
Let $0 < \alpha < \beta<\gamma < 2 \pi\,.$ The four points
$1, e^{i \alpha},e^{i \beta},e^{i \gamma} $ determine a hyperbolic
quadrilateral, whose vertices these points are and whose sides
are orthogonal arcs terminating at these points \cite{lv,hrv1}. We 
consider the problem of finding the modulus (or a lower bound for it) of the 
family $\Gamma$ of curves within the quadrilateral joining
the opposite orthogonal arcs $ (e^{i \alpha},e^{i \beta} )$ and
$ (e^{i \gamma},1 )$ within the quadrilateral \cite{lv}\,. It is
easy to see that we can find a M\"obius transformation $h$ of
$\mathbb{D}$ onto
$\IH$ such that $h(1) =1$, $h(e^{i \alpha}) =t,$ $h(e^{i \beta})= 
-t$, $h(e^{i \gamma}) = -1 $ for some $t>1\,.$  The number $t$ can
be found by setting the absolute ratios $|1, e^{i \alpha},e^{i 
\beta},e^{i \gamma}|$ and $|1, t, -t,1|$ equal, and solving the resulting
quadratic equation for $t$ because M\"obius transformations preserve
absolute ratios. The image quadrilateral has four semicircles as its
sides, the diameters of these are $[-1,1], [1,t],[-t,t],[-t, -1]\,$
and the family $h(\Gamma)$ has a subfamily $\Delta$ consisting of 
radial
segments
$$
[e^{i \phi}, te^{i \phi}], \quad \phi\in(\theta, \pi -\theta),\quad
\sin \theta =\frac{t-1}{t+1}\,.
$$
Obviously, for $\theta =0$ we obtain an upper bound. Therefore
$$
\frac{\pi}{\log t} \ge \M(h(\Gamma))\ge \M(\Delta)= \frac{\pi- 2 \theta}{ \log t}  \,.
$$

\subsubsection{Numerical Experiments}\label{sec:hyperbolicexperiments}
\begin{figure}
\centering
\subfloat[Case 1: $s=\pi/4$.]{\includegraphics[height=2in]{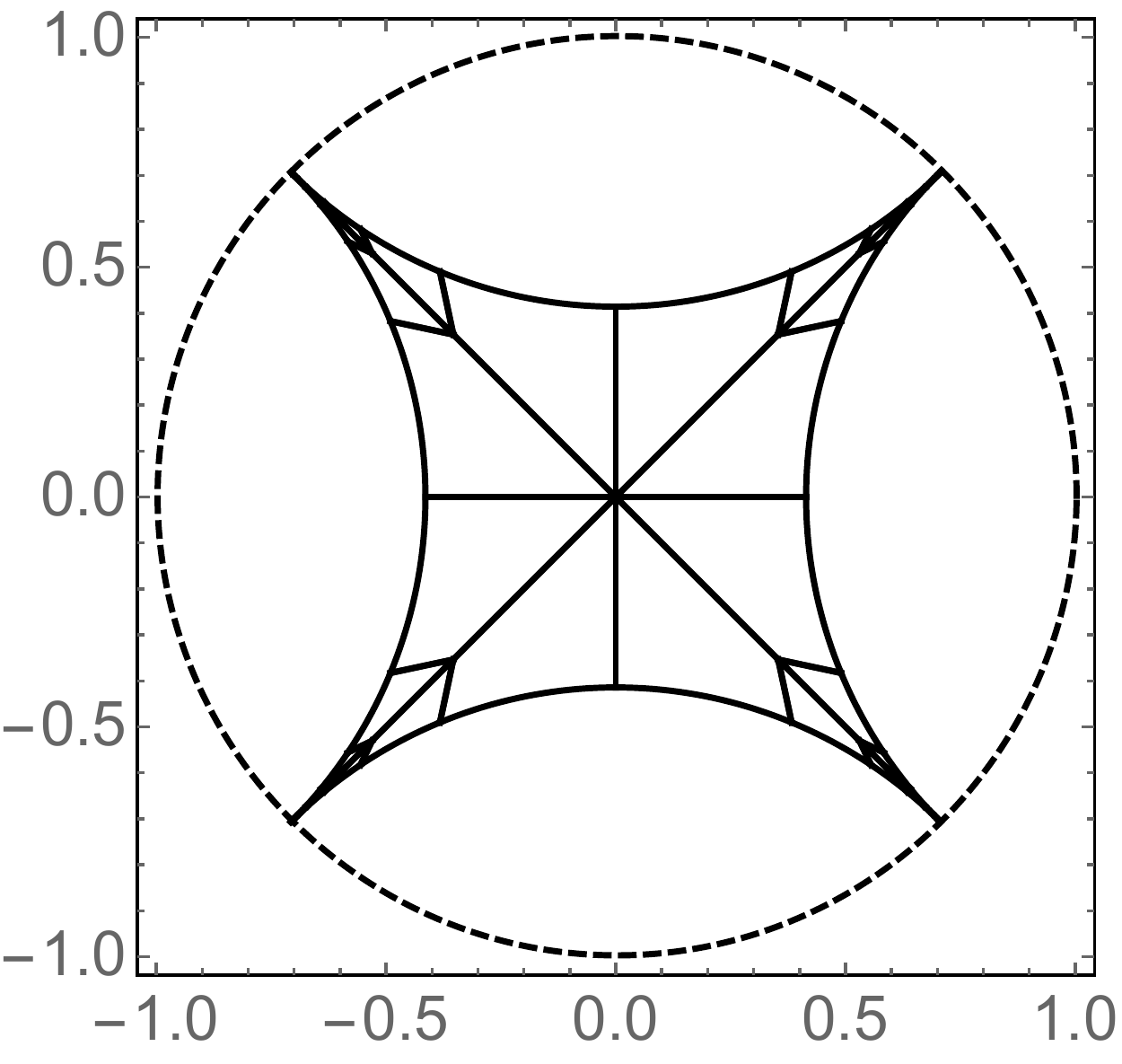}}\quad
\subfloat[Case 2: $s=3\pi/8$.]{\includegraphics[height=2in]{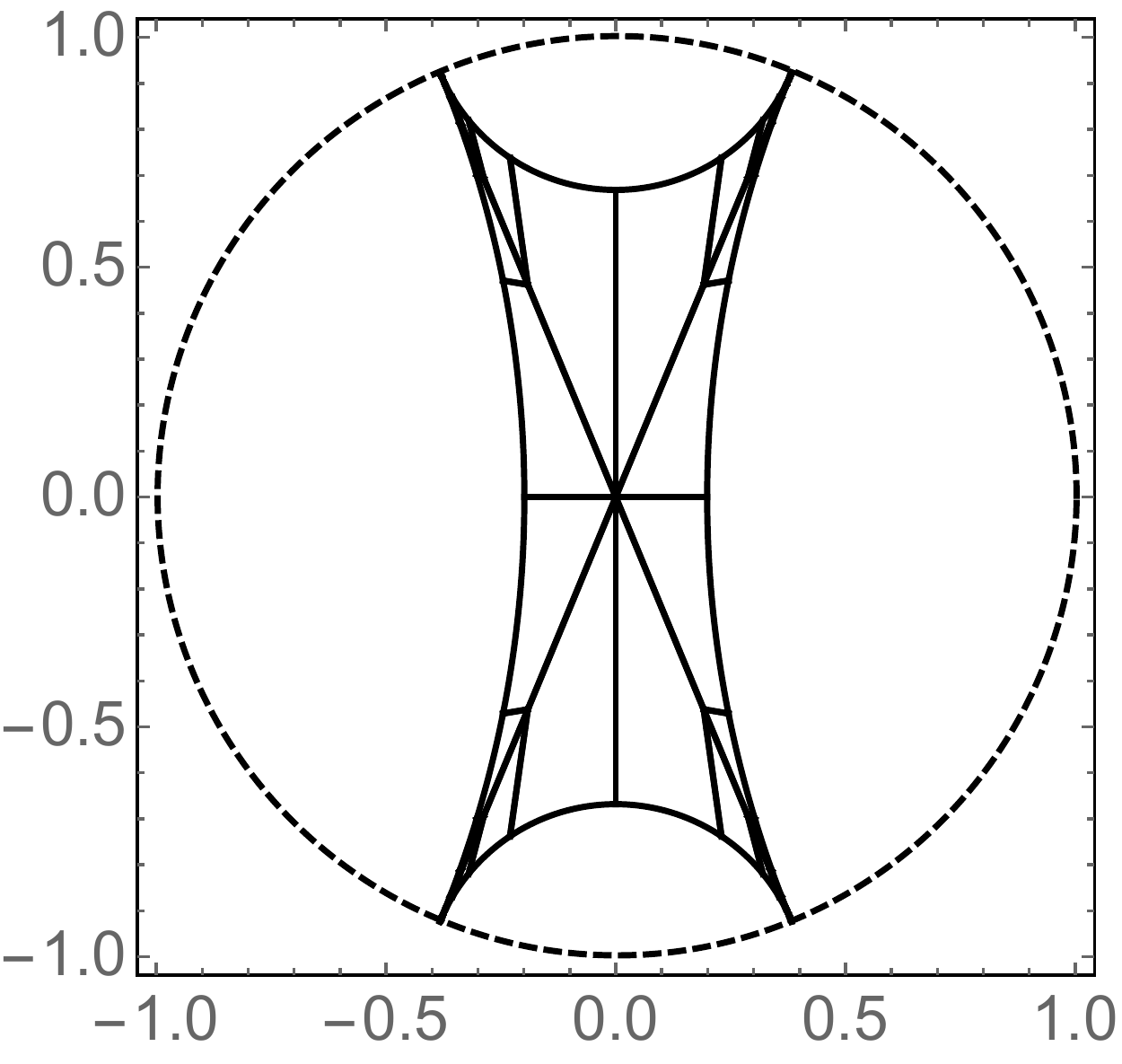}}
\caption{Hyperbolic rectangles: Refined meshes.}
\label{fig:hypdomain}
\end{figure}

\begin{figure}
\centering
\subfloat[Case 1.]{\includegraphics[height=2in]{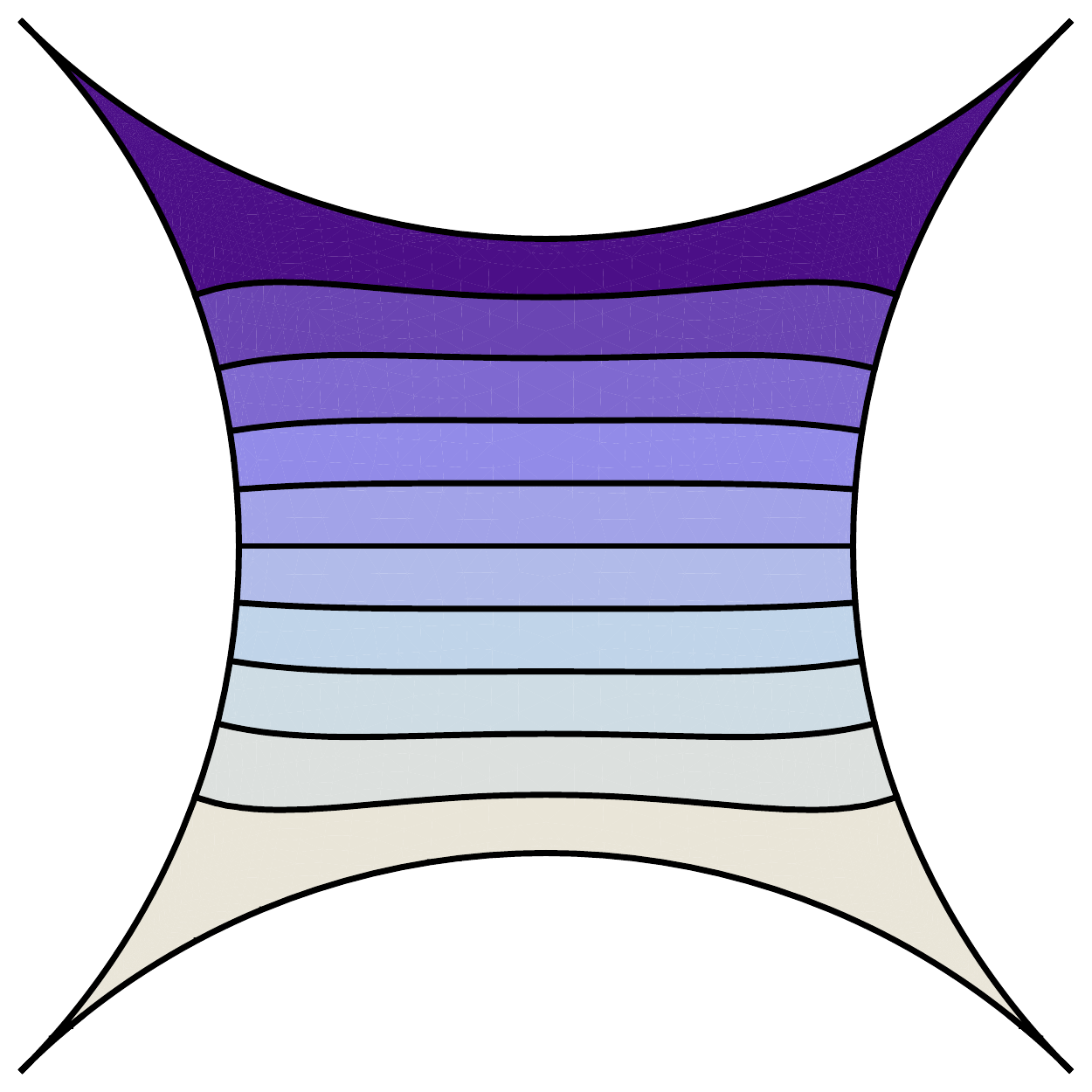}}\quad
\subfloat[Case 2.]{\includegraphics[height=2in]{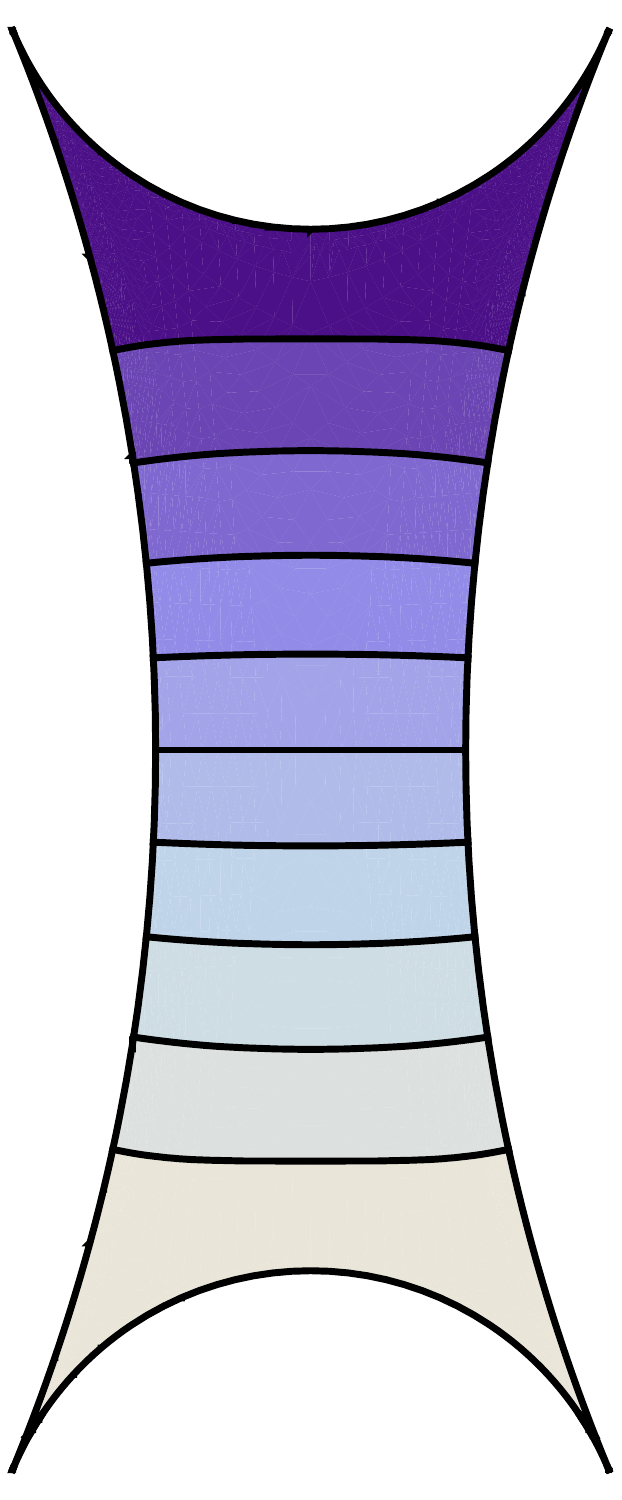}}\quad
\subfloat[Case 2 (Conjugate).]{\includegraphics[height=2in]{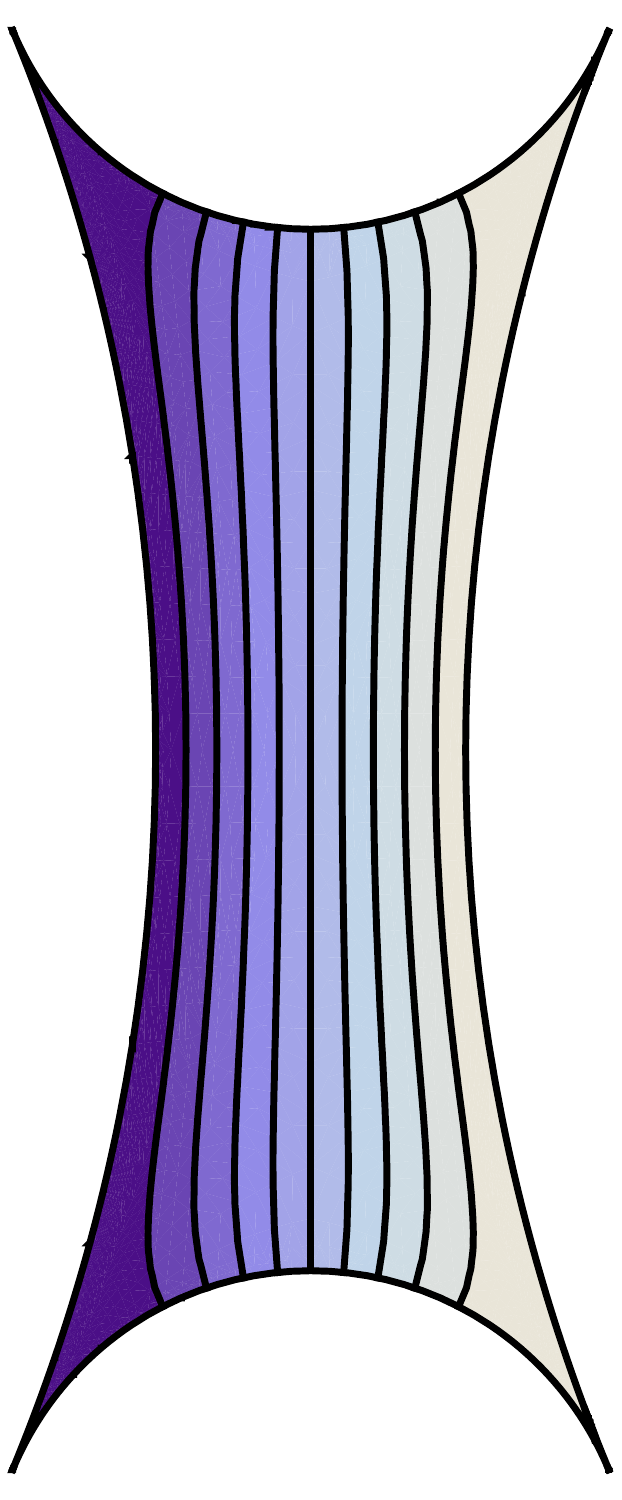}}
\caption{Hyperbolic rectangles: Potential functions.}
\label{fig:hyppotential}
\end{figure}

\begin{figure}
\centering
\subfloat[Case 1.]{\includegraphics[width=0.45\textwidth]{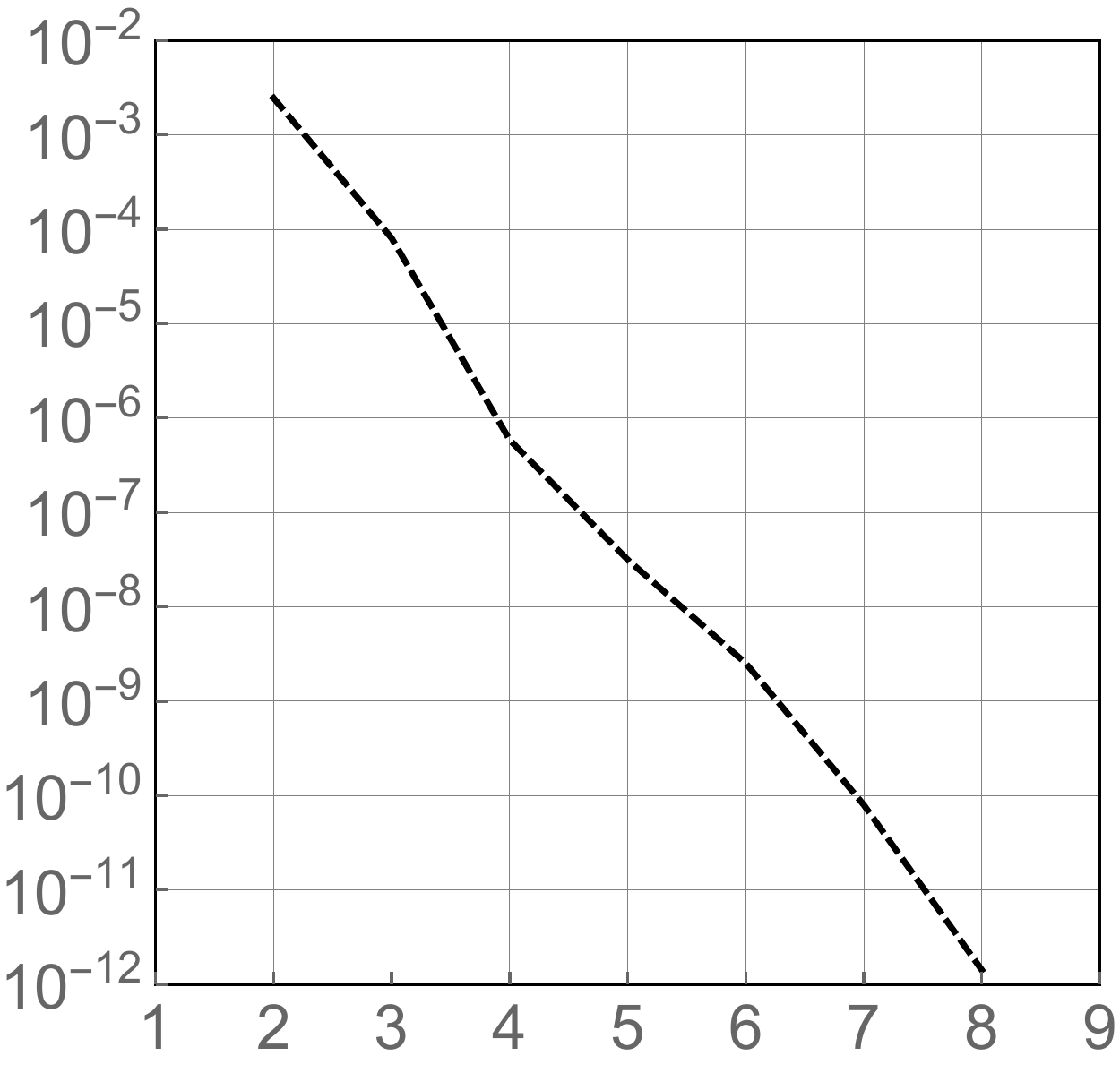}}\quad
\subfloat[Case 2.]{\includegraphics[width=0.45\textwidth]{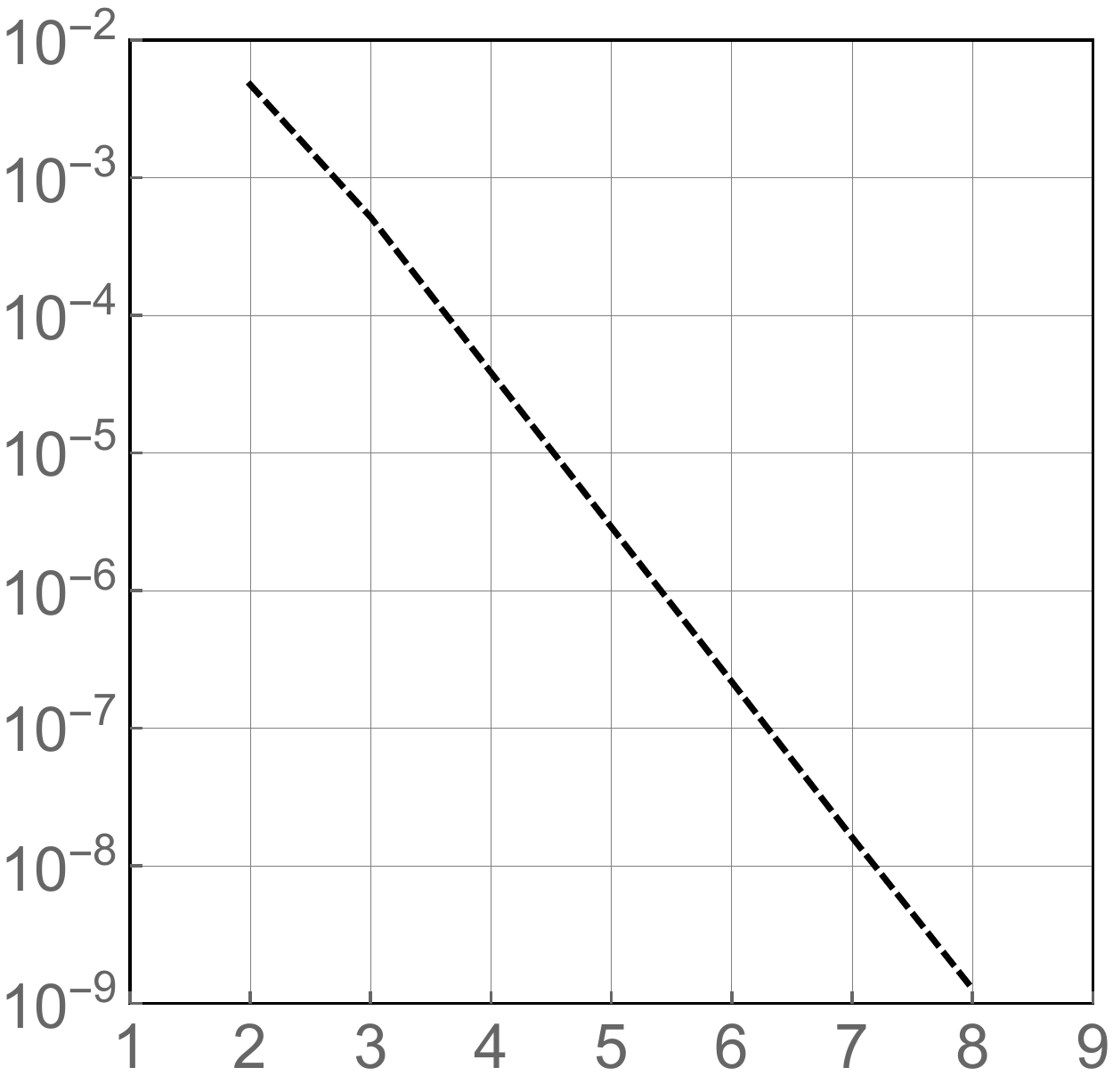}}
\caption{Hyperbolic rectangles: Reciprocal errors; $\log$-plot: Error vs $p$.}
\label{fig:hyprecip}
\end{figure}

Similarly as before, the examples of this section are outlined in Table~\ref{tbl:komori}
and Figures~\ref{fig:hypdomain}, \ref{fig:hyppotential}. The meshes are refined
in exactly same fashion so that any differences in convergence stem only from the
difference in the geometric scaling. As shown in Figure~\ref{fig:hyprecip}
the convergence in the reciprocal error is exponential, but with a better rate
for the symmetric case. Moreover, for the symmetric
domain both error estimates coincide.

\begin{figure}
\centering
\subfloat[Case 1.]{\includegraphics[width=0.3\textwidth]{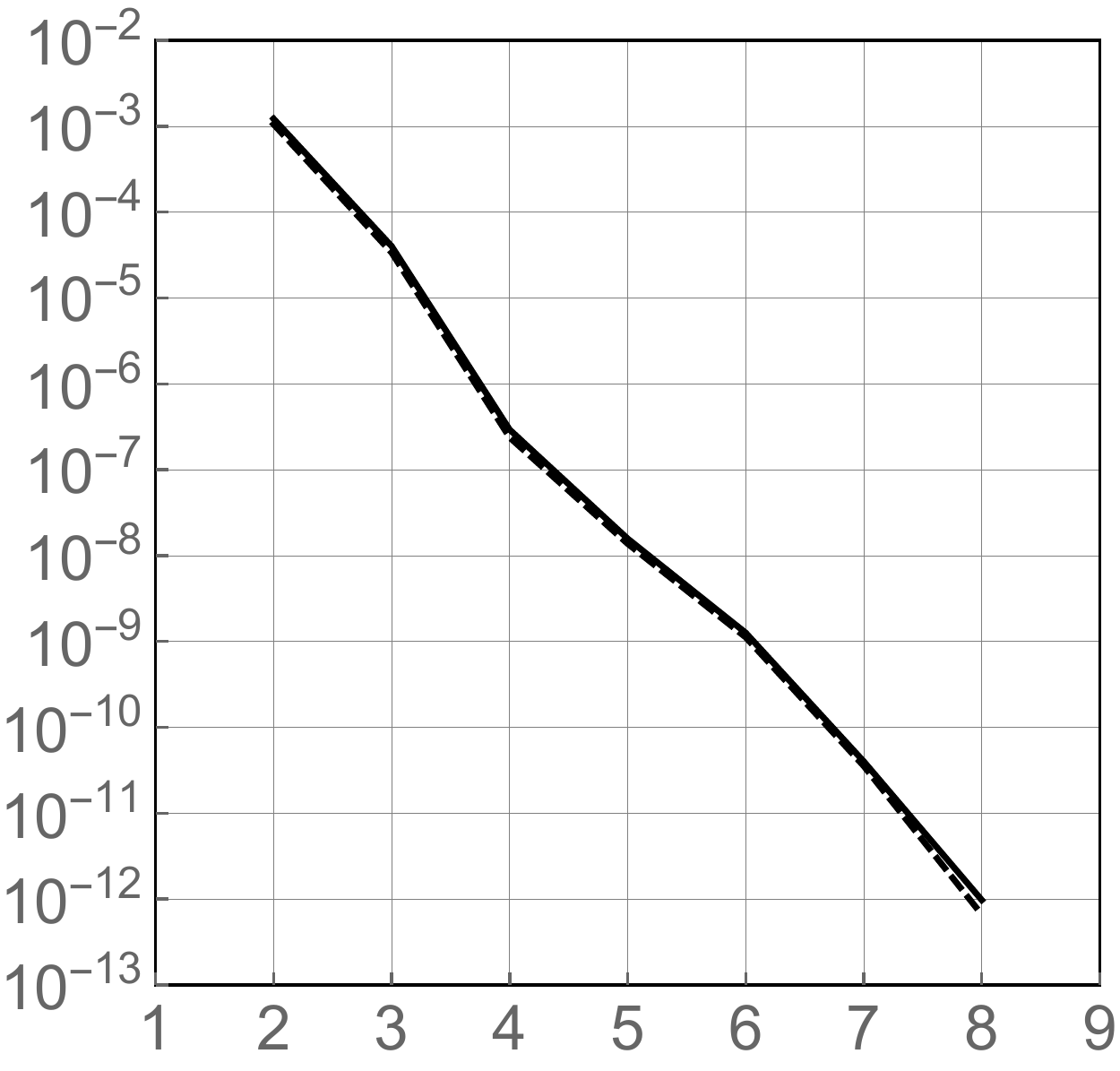}}\quad
\subfloat[Case 2.]{\includegraphics[width=0.3\textwidth]{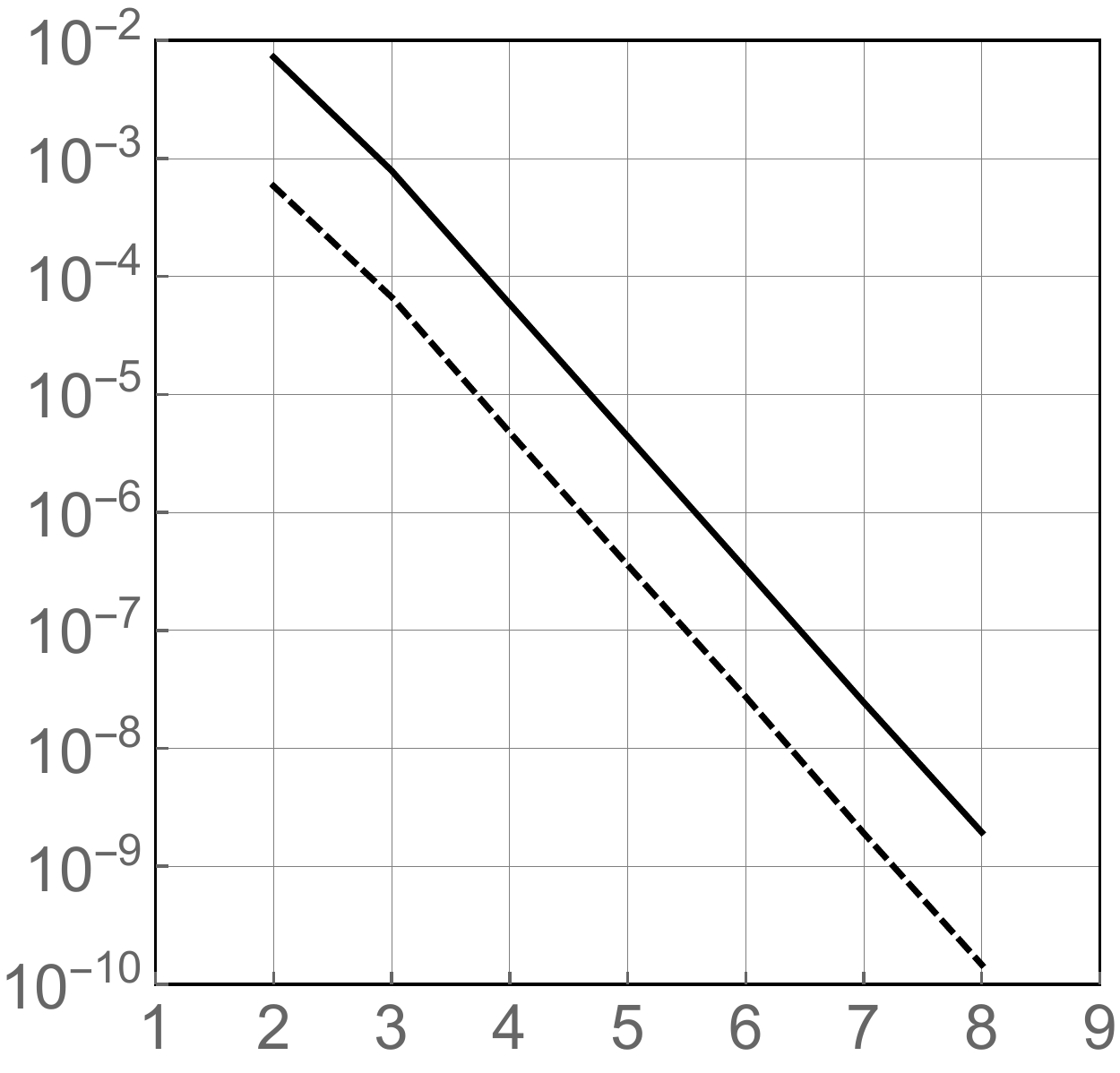}}\quad
\subfloat[Case 2 (Conjugate).]{\includegraphics[width=0.3\textwidth]{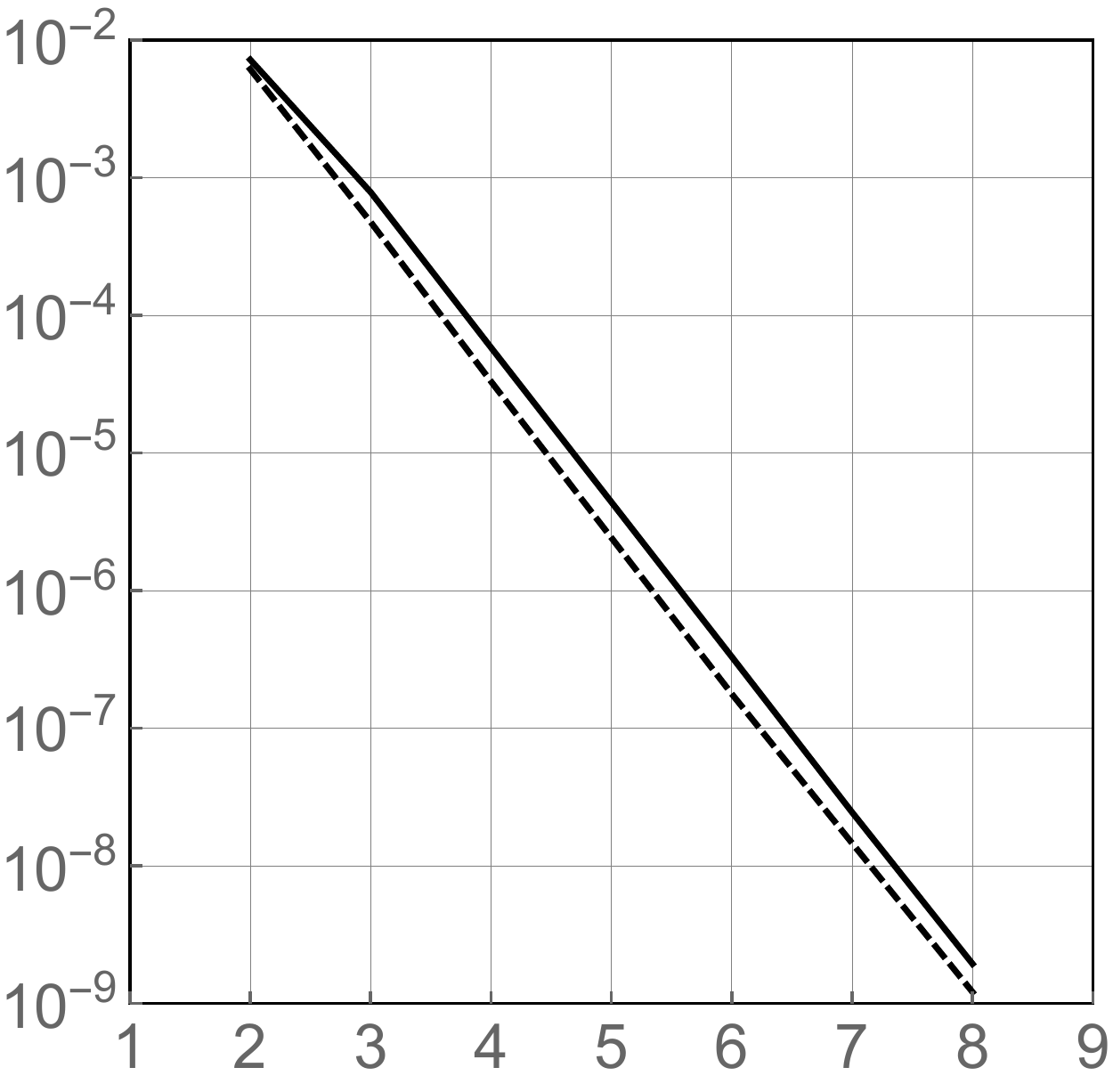}}
\caption{Hyperbolic rectangles: Estimated errors; $\log$-plot: Error vs $p$; Solid line = Reciprocal estimate, Dashed line = Auxiliary space estimate.}
\label{fig:hyperror}
\end{figure}

\begin{table}[ht]
\caption{
  Tests on hyperbolic quadrilaterals. The errors are given as {$|\lceil\log_{10}|\mathrm{error}|\rceil|$}.
}\label{tbl:komori}
\renewcommand\arraystretch{1}
\noindent

\begin{center}\footnotesize
\begin{tabular}{|l|l|l|l|l|l|}
\hline
Case & Method & Errors & Sizes & $M(Q_s)$\\
\hline
1 & $hp$, $p=12$ & 12 & 10225 & 1\\
2 & $hp$, $p=16$ & 11 &17985 & 3.037469188986459\\
\hline
\end{tabular}
\end{center}
\end{table}

\subsection{Other Quadrilaterals with Cusps}
\label{sec:quadrilateralswithcusps}

Next we consider quadrilaterals  $Q=Q(D; z_1,z_2,z_3,z_4)$ where neither of the components of $\C_\infty \setminus \partial D$ is bounded. The modulus of the following quadrilateral has been obtained by W.~Bergweiler and A.~Eremenko \cite{be}, who studied this question in connection to an extremal problem of geometric function theory introduced by A.A.~Goldberg in 1973.

\subsubsection{Example I}
Consider the strip the closed unit disk is removed:
\[
D_1= \{ z : -3<\re z<1\} \setminus \overline{\mathbb{D}}.
\]
Let the four vertices $z_j$, $j=1,2,3,4$ on the boundary of $D$ be $1,\infty,\infty,1$, in counter-clockwise order. Then, all the angles at vertices are equal to $0$.

First we map the domain in question to a bounded domain so that the line $\{z : \re z = 1\}$ maps to unit circle, and the real axis remains fixed. After the M\"obius transformation we may assume that we are computing in the unit disk $\mathbb{D}$. For convenience, consider disks $\ID(1-t,t)$ and $\ID(-1+s,s)$ internally tangent to the unit circle at the points $1$ and $1$, respectively, with $s,t \in (0, 1/3)$. The corner points of the quadrilateral are $1,-1,-1,1$, with a zero angle at each of the corners. We denote the respective radii of the disks by $s$ and $t$. 

We have computed numerically the modulus of the family of curves joing the
two disks within the domain $D_2 = \ID \setminus (\ID(1-t,t) \cup \ID(-1+s,s))$.  It is the reciprocal of the  modulus of the family of curves joing the upper semicircle with the lower semicircle within the same domain. The results are summarized in Table 

An estimate for the case $s=t=\sqrt{2} -1 \approx 0.41421$ is obtained by Bergweiler and Eremenko, 
with numerical values that agree with our results up to 6 significant digits, in \cite{be}. The conformal modulus in this case is approximately $2.78234$ (AFEM), or $2.7823418086$ ($hp$-FEM, with error number = 10, $p=21$). Our results agree with the result of \cite{be}.

\subsubsection{Example II}
Consider the domain $D$ (a hexagon) in the upper half-plane obtained from the half-strip
\[
\{z=x+iy: 0<x<1,\, 0<y\},
\]
by removing two half-disks
\[
C_1=\overline{\mathbb{D}(7/24,1/24)},\quad C_2=\overline{\mathbb{D}(5/12,1/12)},
\]
where $\mathbb{D}(z,r)$ denotes the disk centered at $z\in\mathbb{C}$ with radius $r>0$.
Note that $C_1\cap \mathbb{R}=[1/4,1/3]$ and $C_2\cap \mathbb{R}=[1/3,1/2]$. We compute the moduli of two quadrilaterals:
\[ 
Q_1=(D;\infty,0,1/2,1), \quad Q_2=(D; 0,1/4,1/2,1).
\]

Again, we first use the M\"obius transformation
\[
z \mapsto \frac{2z-1}{2z+1}
\]
to map the domain $D$ in question to a bounded domain. Then the boundary points of $Q_1$ map onto
the points $1, -1, 0, 1/3$, respectively. For $Q_2$, the boundary points are mapped onto the points $-1,-1/3,0,1/3$. Quadrilaterals $Q_1$ and $Q_2$, after the M\"obius transformation, and the corresponding potential functions are illustrated in Figure \ref{eremenko-fig2}.

\subsubsection{Numerical Experiments}\label{sec:eremenkoexperiments}
\begin{figure}
\centering
\subfloat[Case 1: $s=t=\sqrt{2} -1$.]{\includegraphics[width=0.45\textwidth]{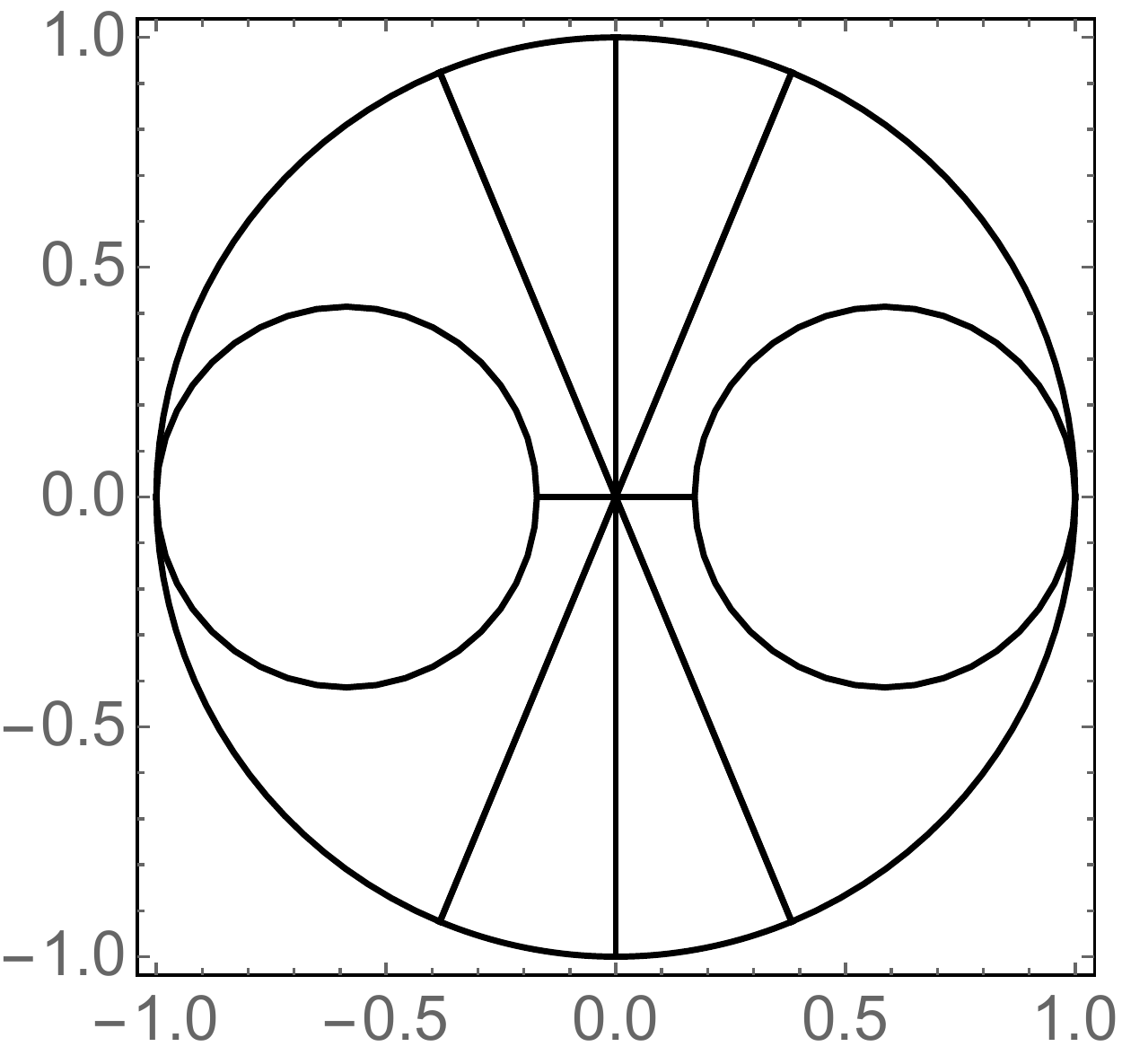}}\quad
\subfloat[Case 2: $s=3/10, t=2/5$.]{\includegraphics[width=0.45\textwidth]{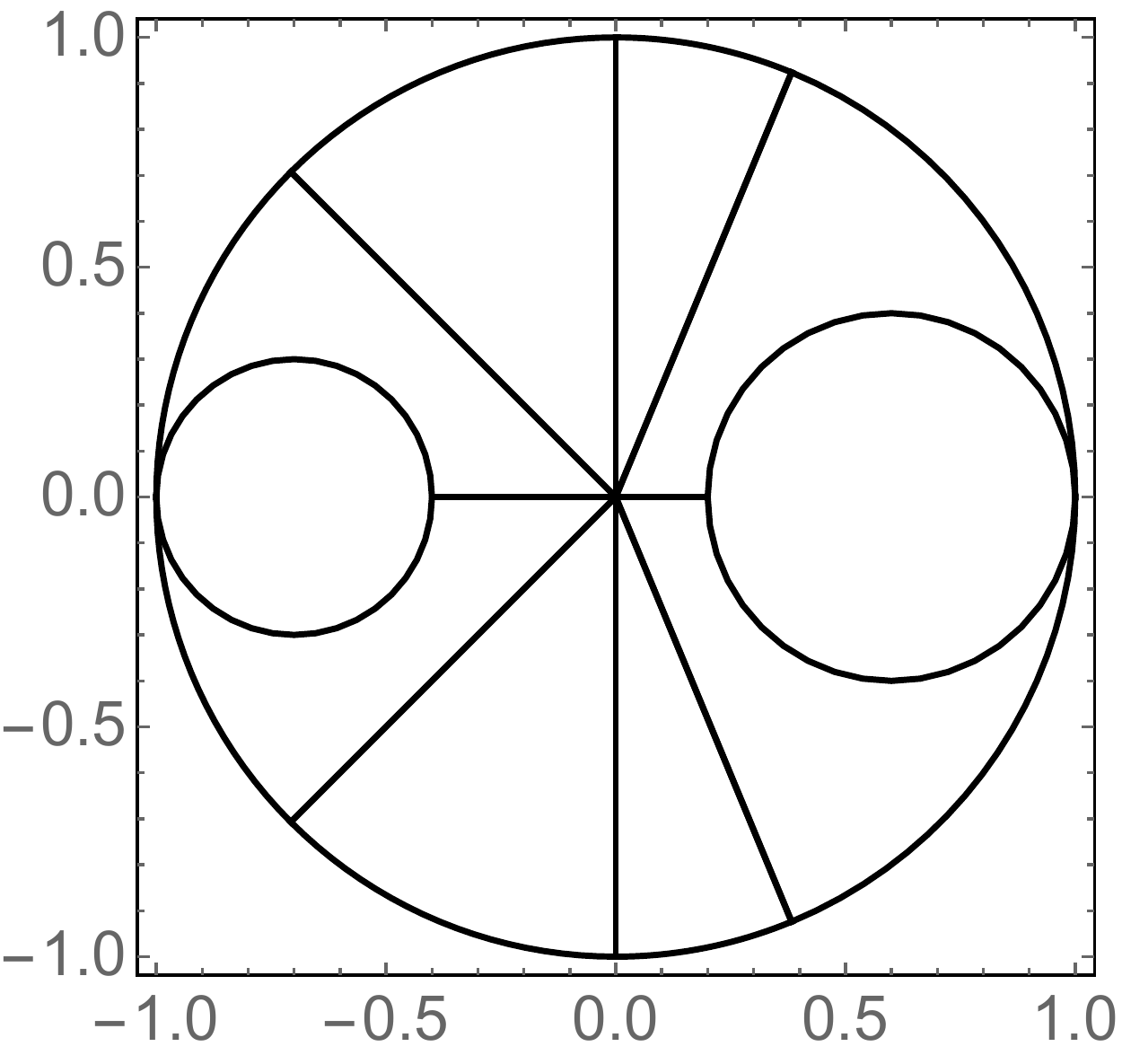}}\\
\subfloat[Case 3 and 4: $s=1/15, t=1/10, u=1/3$.]{\includegraphics[width=0.45\textwidth]{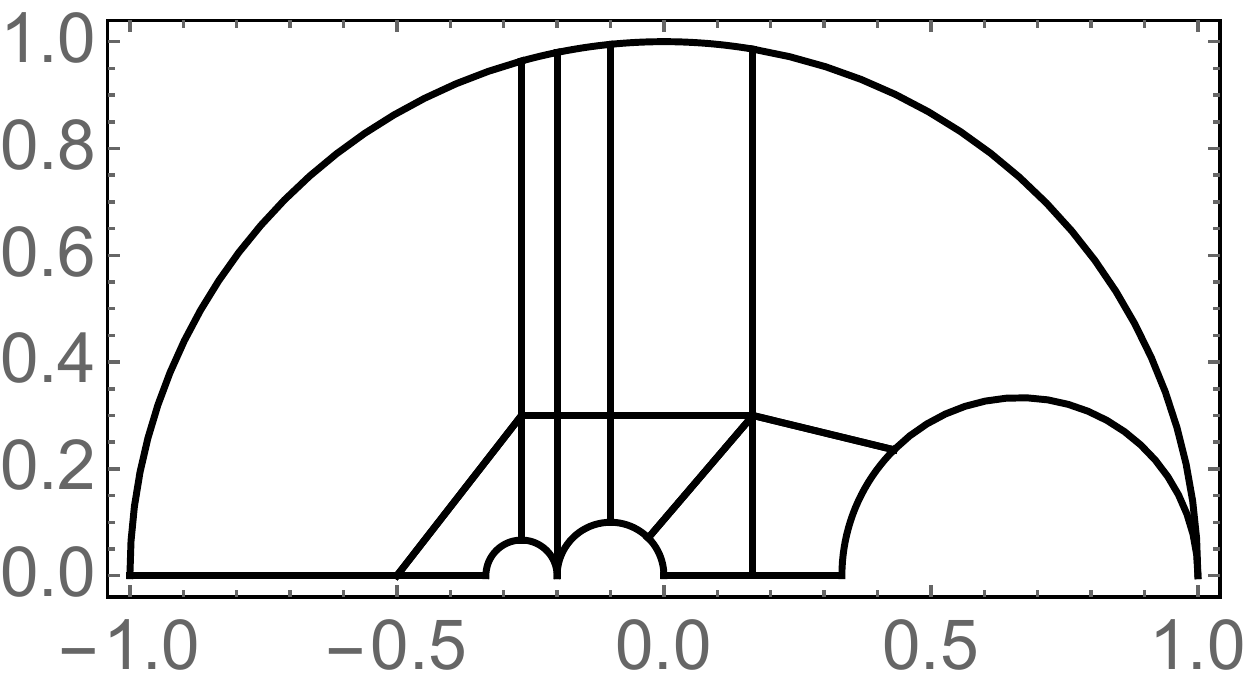}}
\caption{Quadrilaterals with Cusps: $p$-type Meshes.}
\label{fig:cuspdomain}
\end{figure}

\begin{figure}
\centering
\subfloat[Case 1.]{\includegraphics[width=0.45\textwidth]{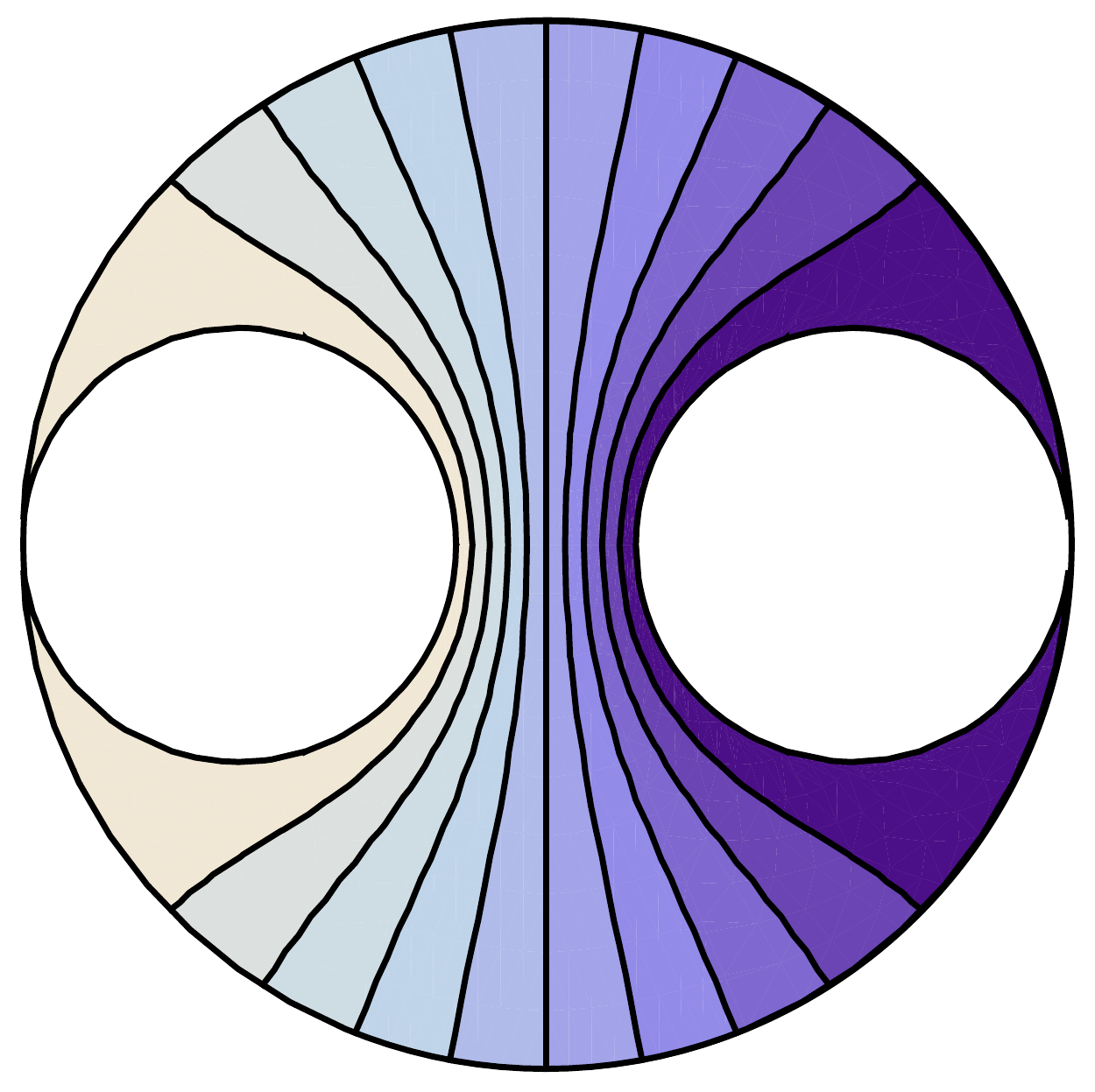}}\quad
\subfloat[Case 2.]{\includegraphics[width=0.45\textwidth]{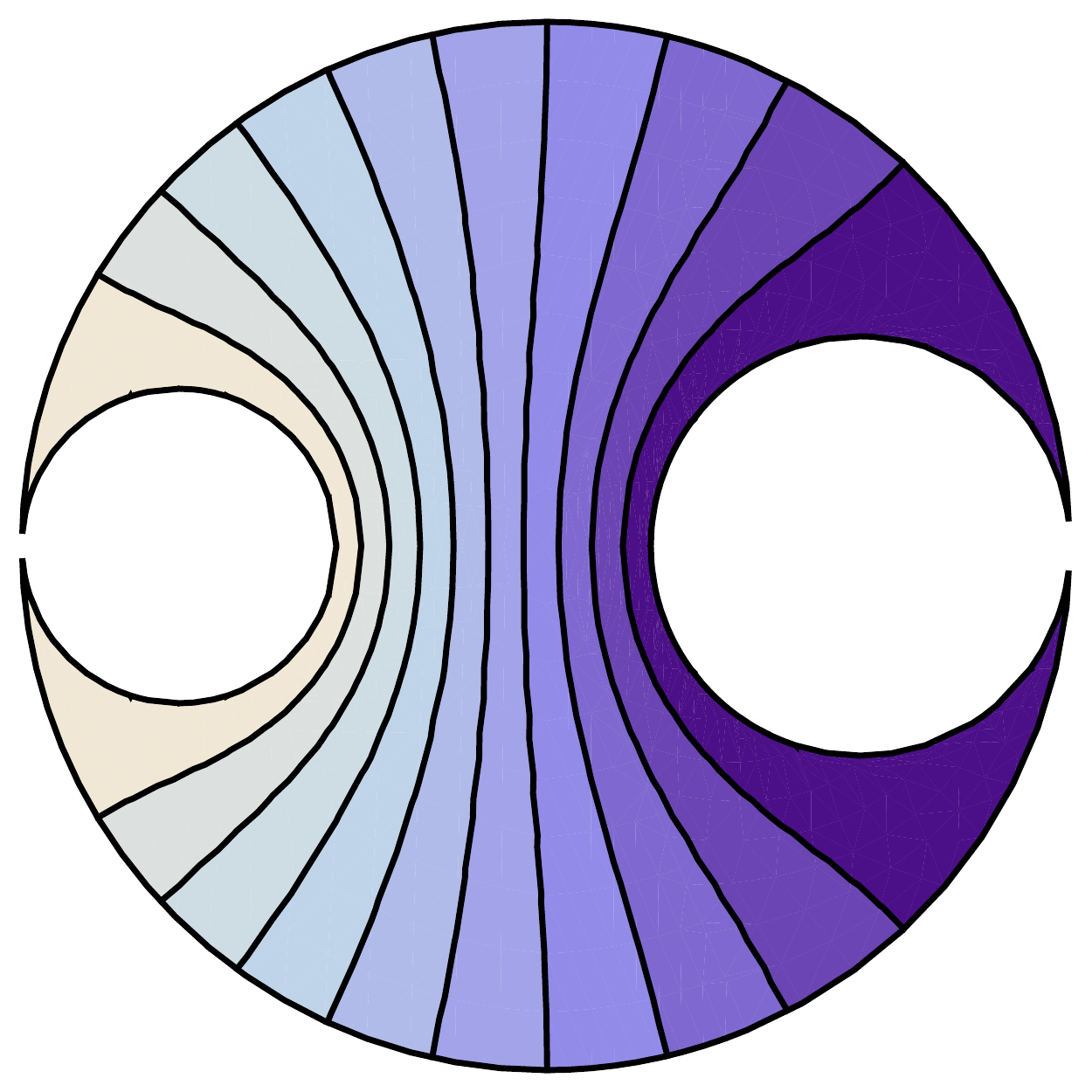}}\\
\subfloat[Case 3.]{\includegraphics[width=0.45\textwidth]{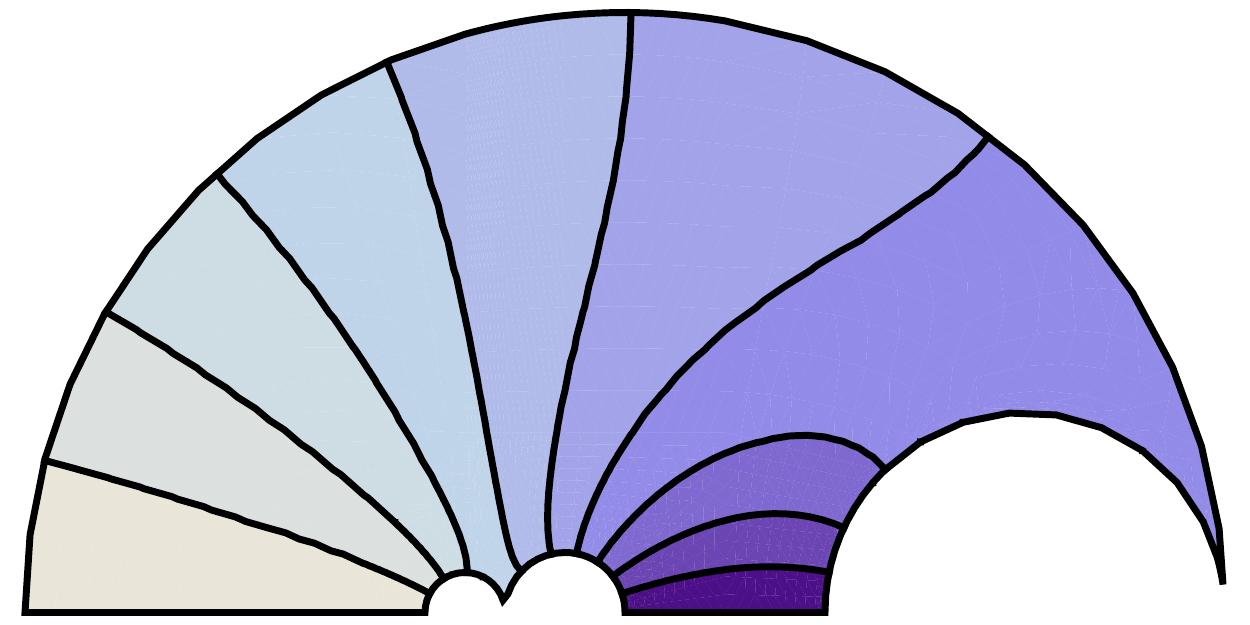}}\quad
\subfloat[Case 4.]{\includegraphics[width=0.45\textwidth]{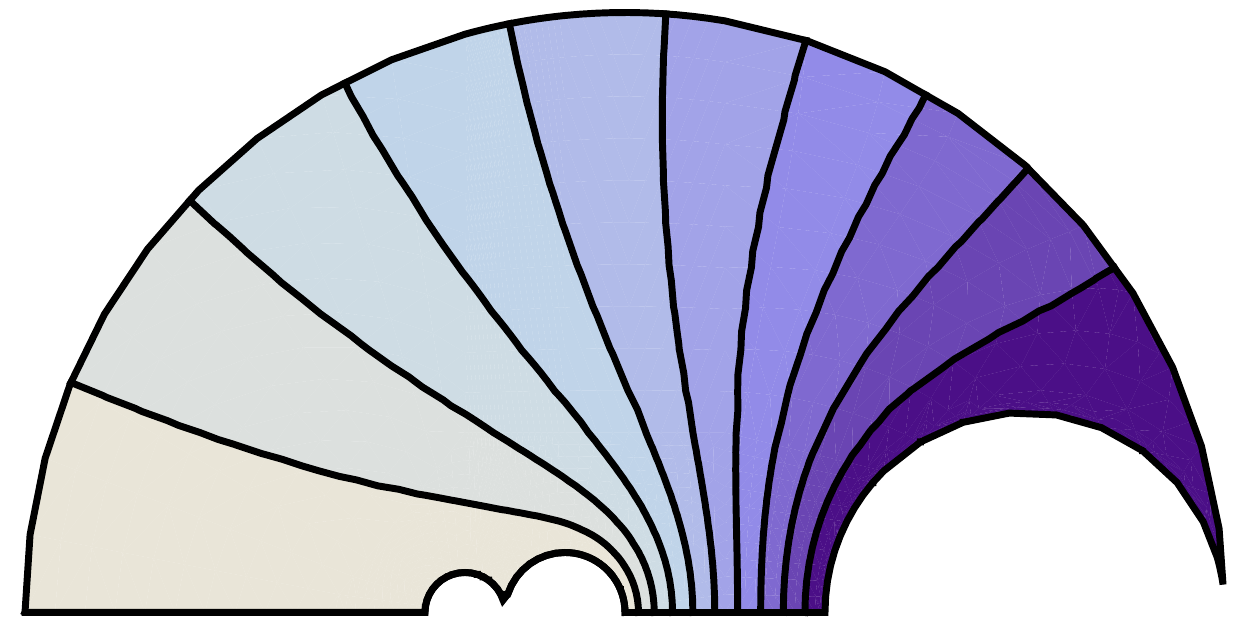}}
\caption{Quadrilaterals with Cusps: Potentials.}\label{eremenko-fig2}
\end{figure}

\begin{figure}
\centering
\subfloat[Reciprocal Error.]{\includegraphics[width=0.45\textwidth]{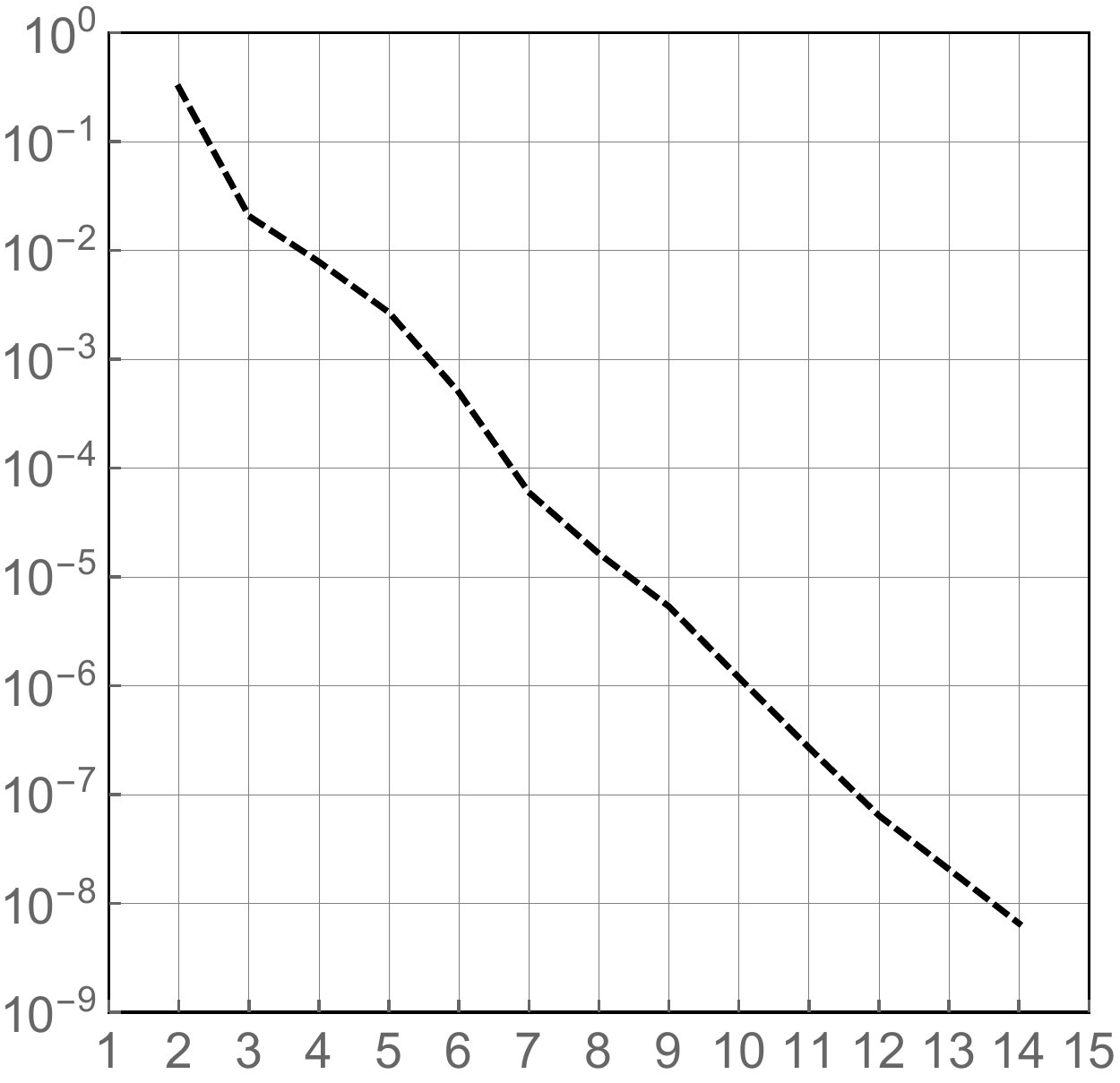}}\quad
\subfloat[Estimated Error.]{\includegraphics[width=0.45\textwidth]{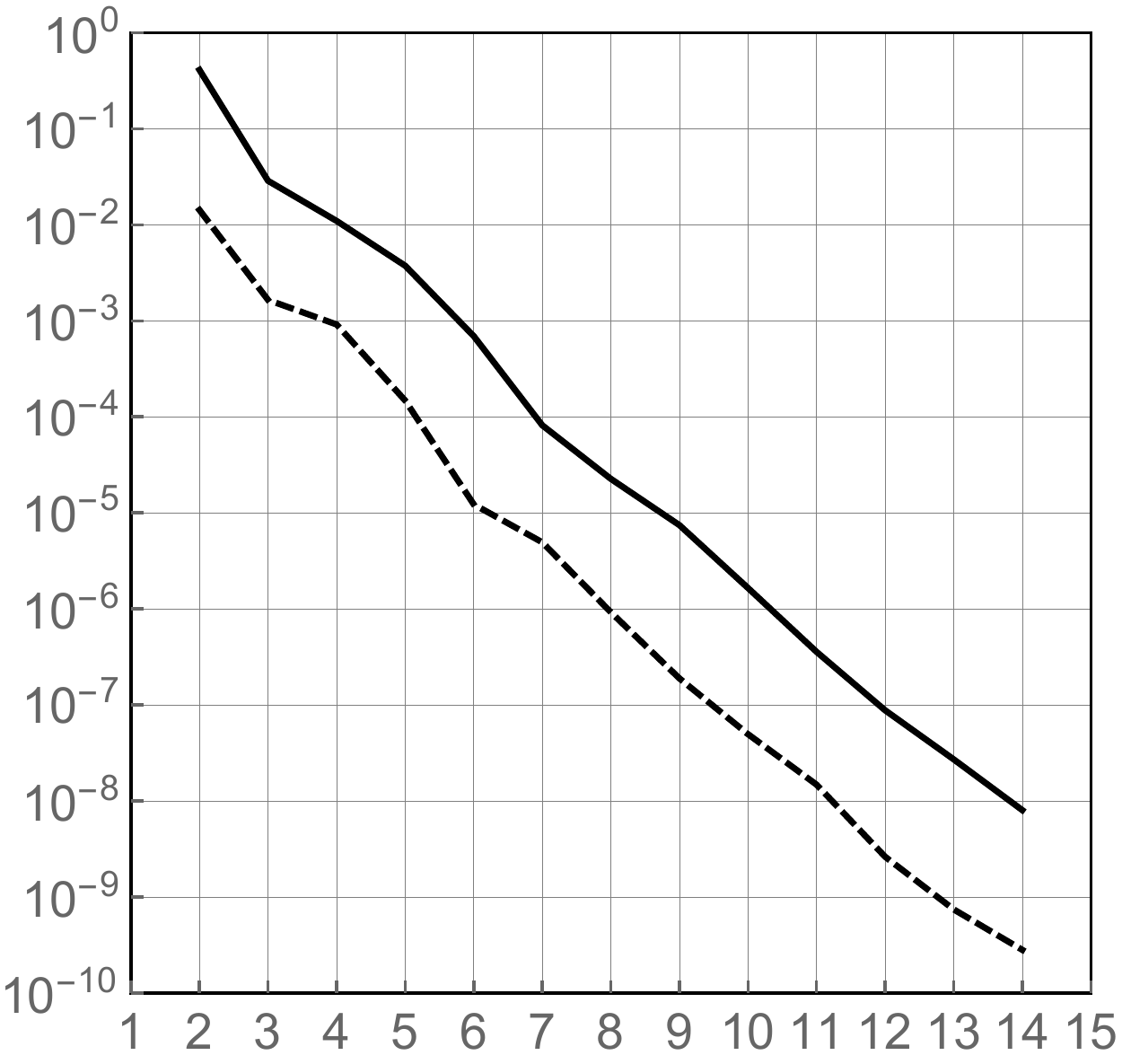}}\\
\subfloat[Estimated Error (Conjugate).]{\includegraphics[width=0.45\textwidth]{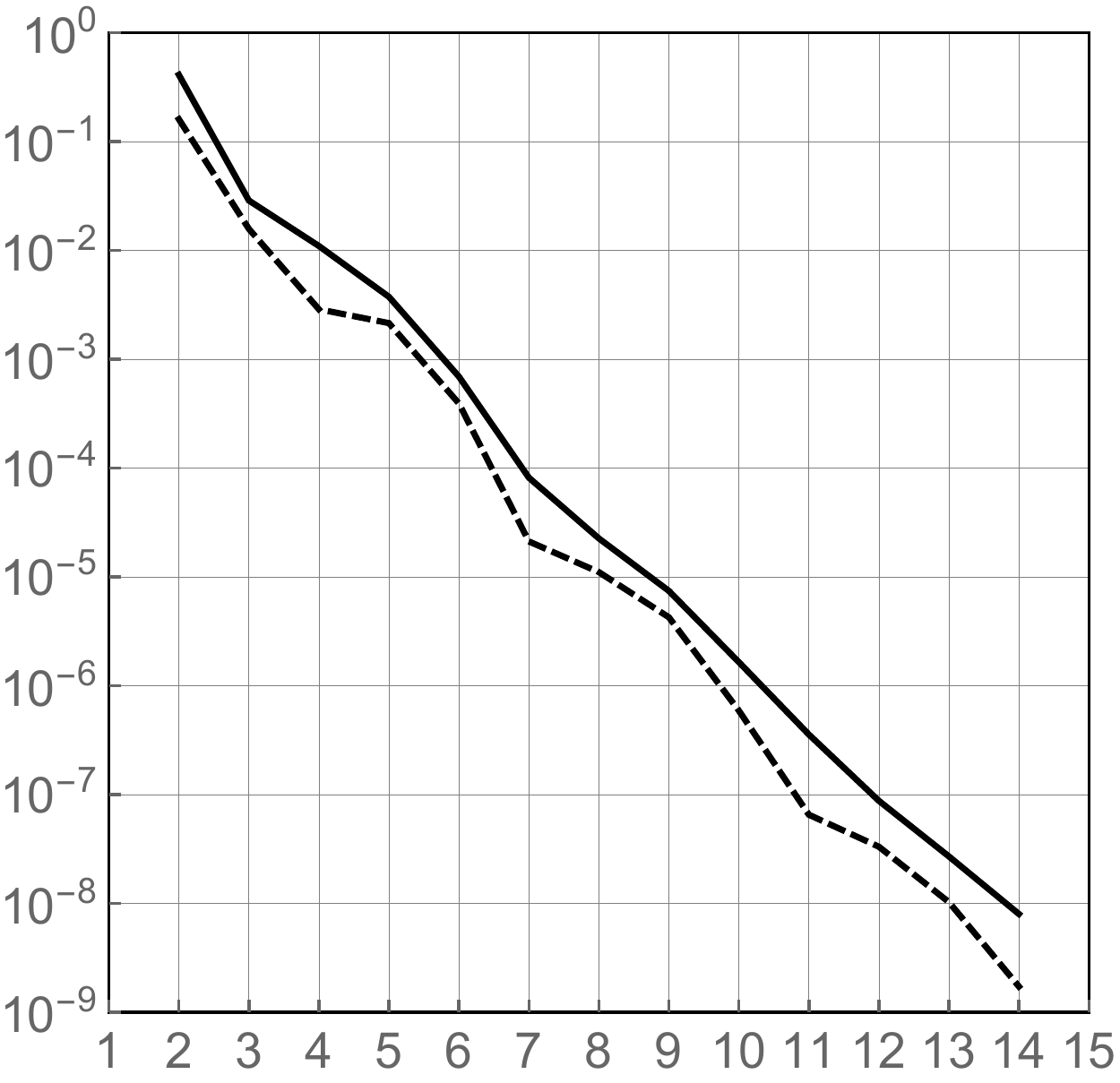}}\quad
\subfloat[Estimated Error: $(e,b)=(0,2)$.]{\includegraphics[width=0.45\textwidth]{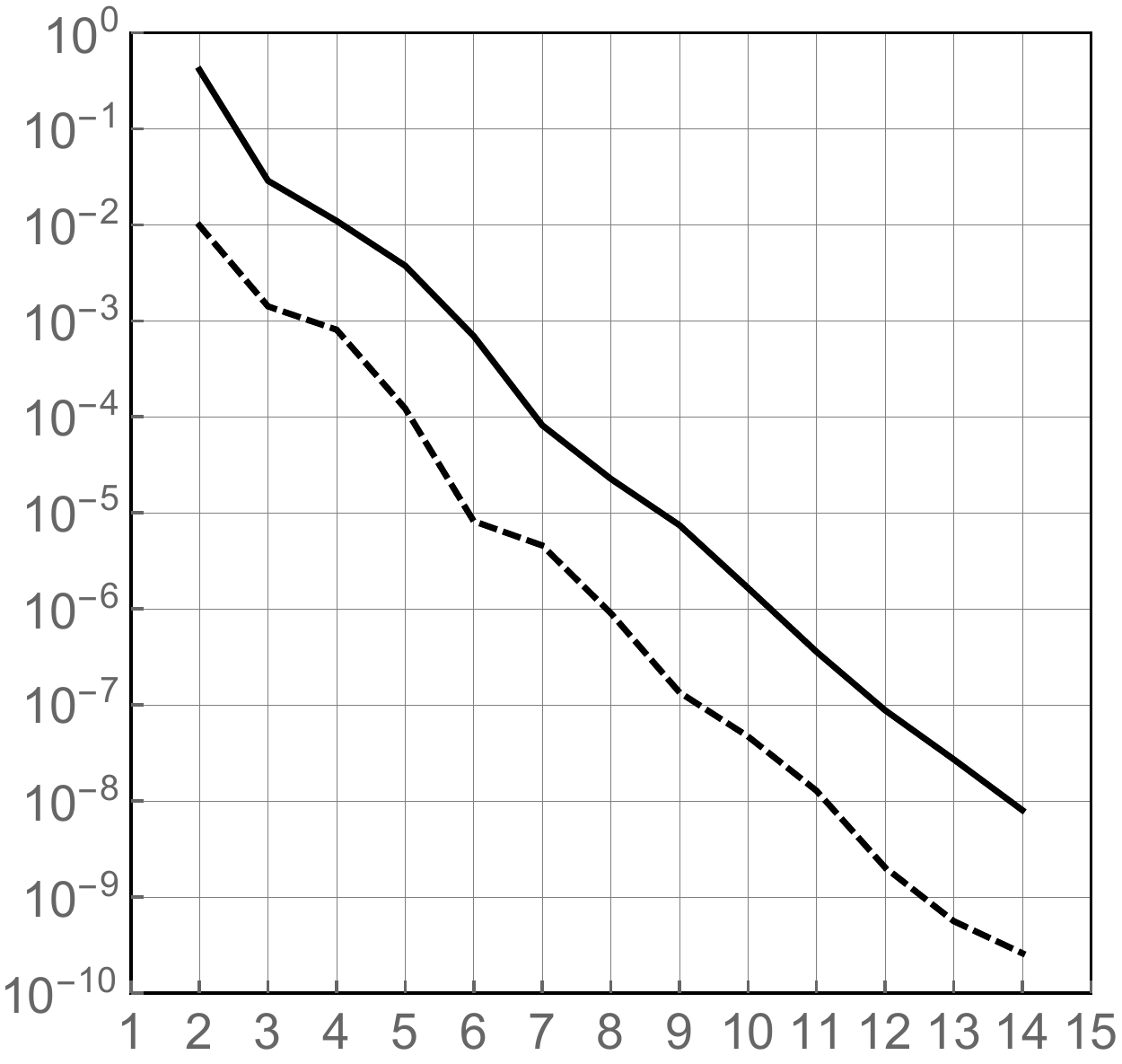}}
\caption{Quadrilaterals with Cusps: Case 1: Reciprocal error, $\log$-plot: Error vs $p$;  Estimated error, $\log$-plot: Error vs $p$; Solid line = Reciprocal estimate, Dashed line = Auxiliary space estimate.}
\label{fig:cuspcase0}
\end{figure}

\begin{figure}
\centering
\subfloat[Case 2.]{\includegraphics[width=0.30\textwidth]{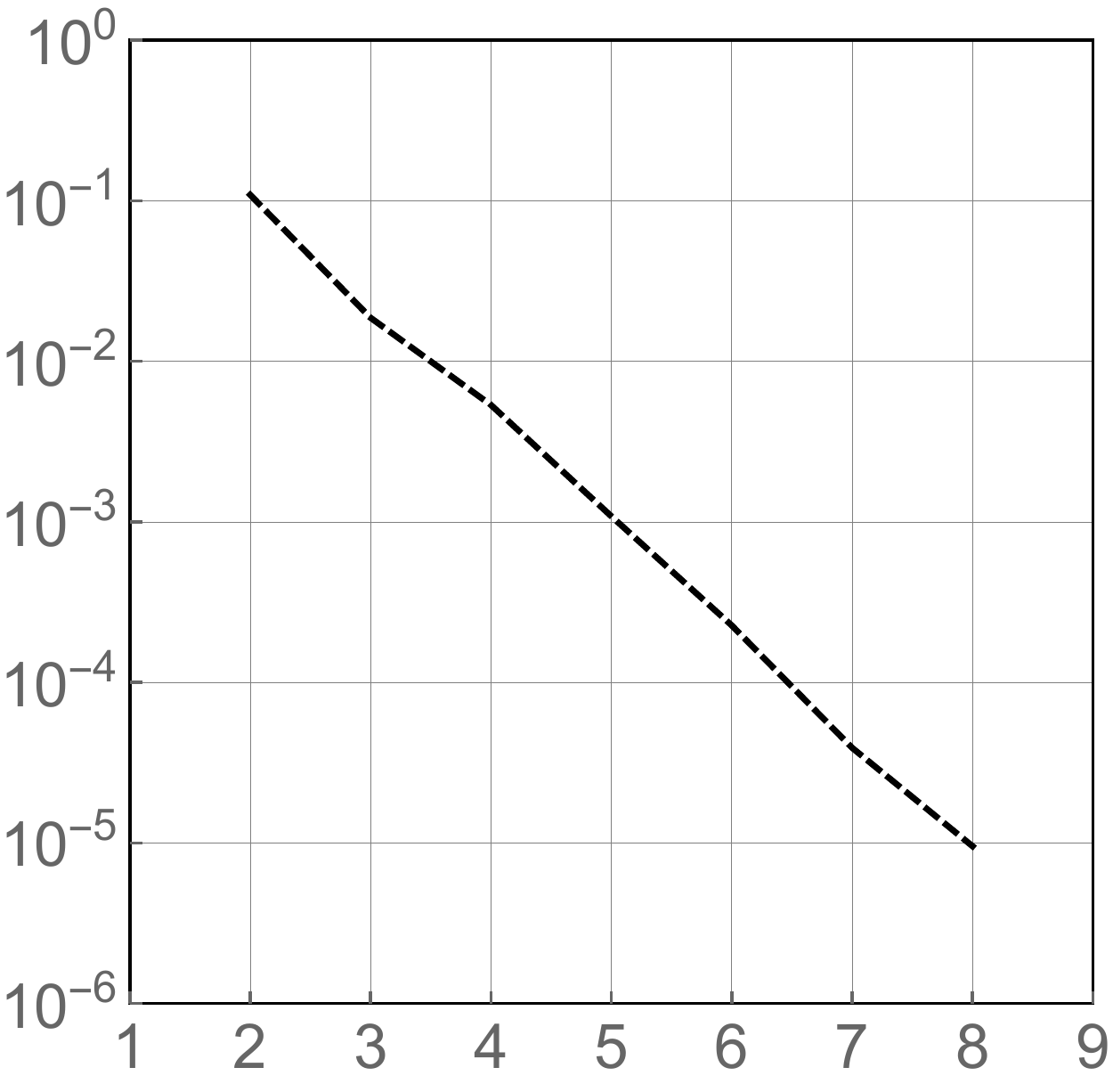}}\quad
\subfloat[Case 3.]{\includegraphics[width=0.30\textwidth]{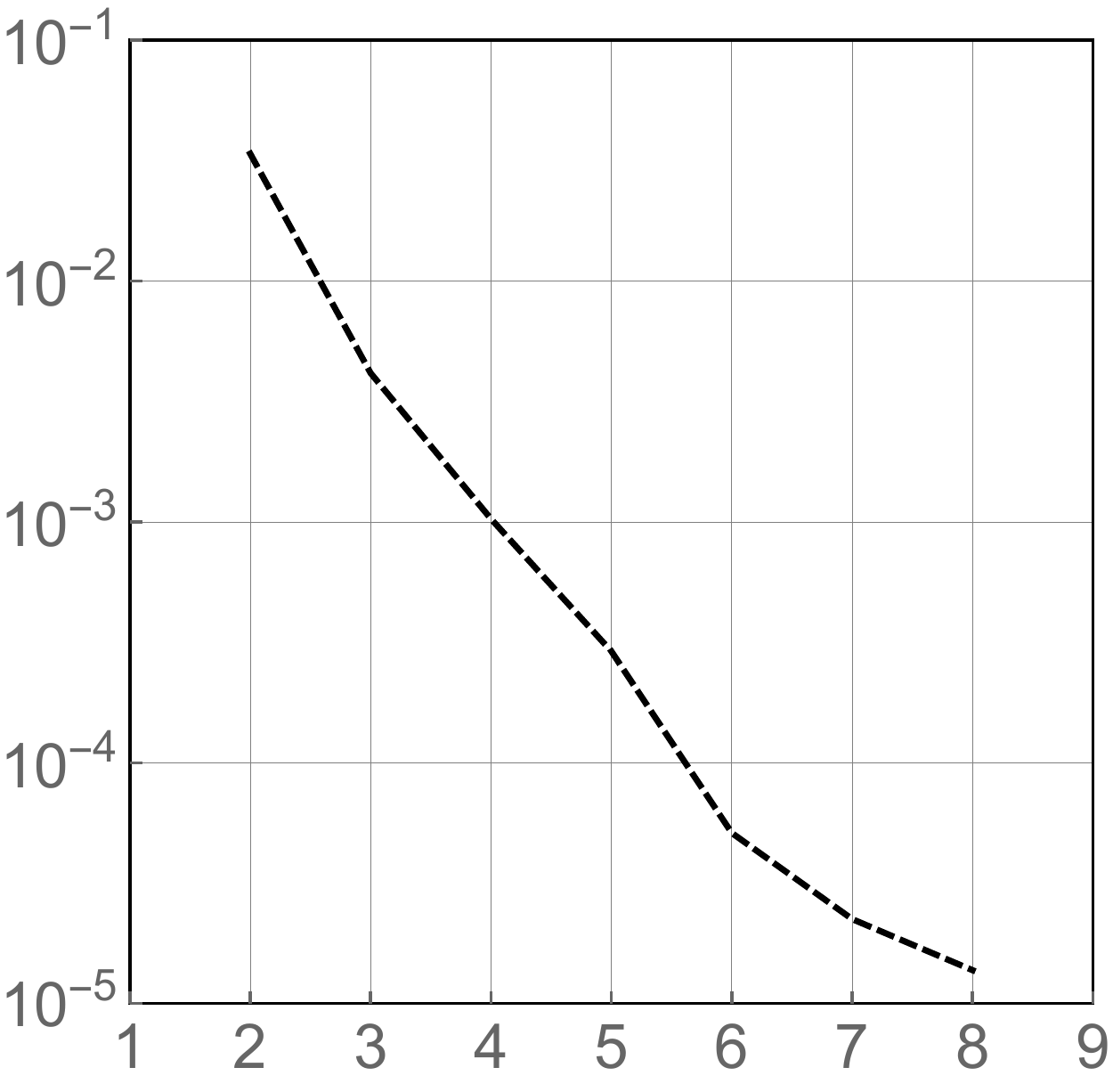}}\quad
\subfloat[Case 4.]{\includegraphics[width=0.30\textwidth]{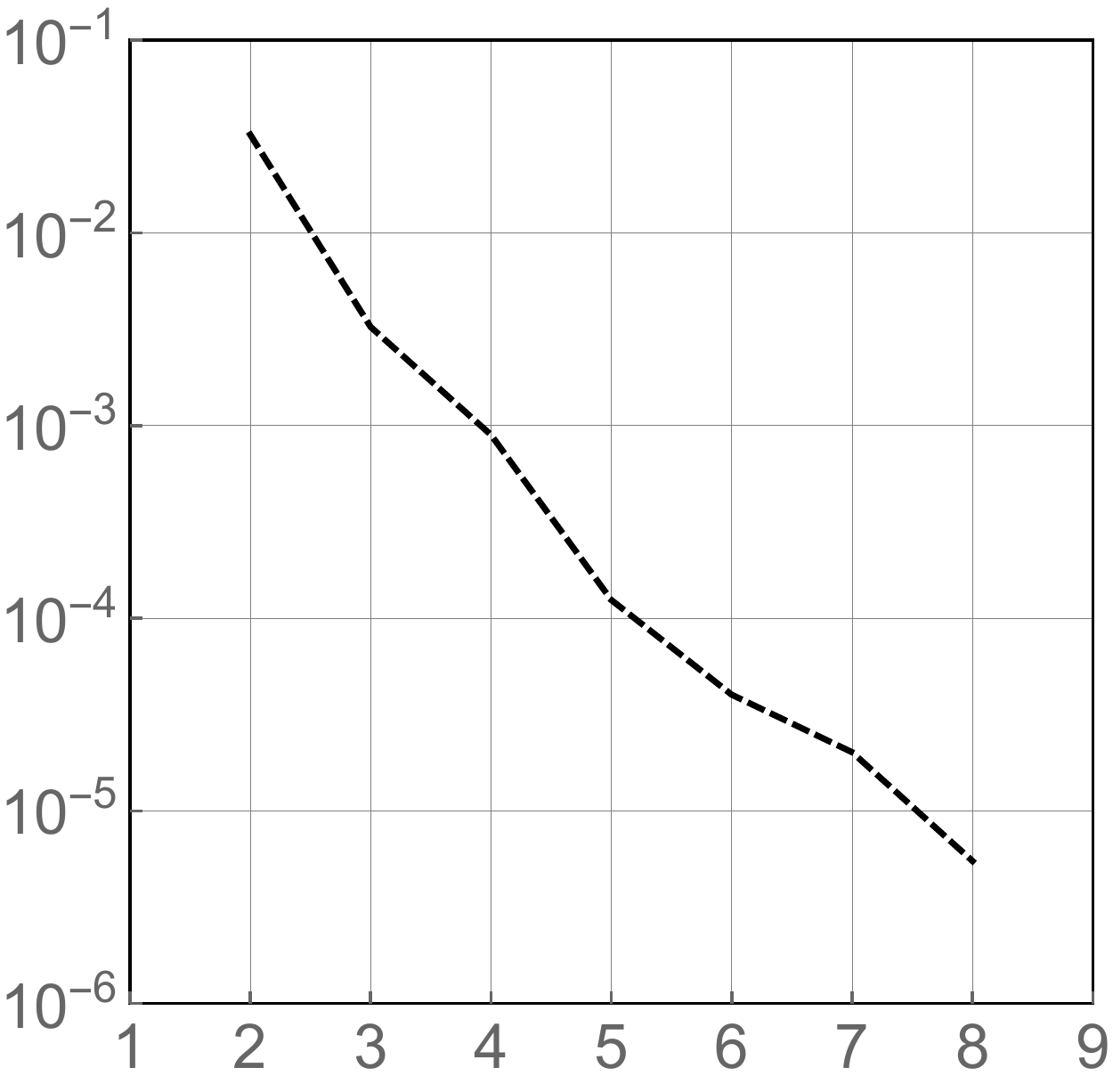}}
\caption{Quadrilaterals with Cusps: Cases 2--4: Reciprocal errors; $\log$-plot: Error vs $p$.}
\label{fig:cuspcase24}
\end{figure}

\begin{figure}
\centering
\subfloat[Case 2.]{\includegraphics[width=0.30\textwidth]{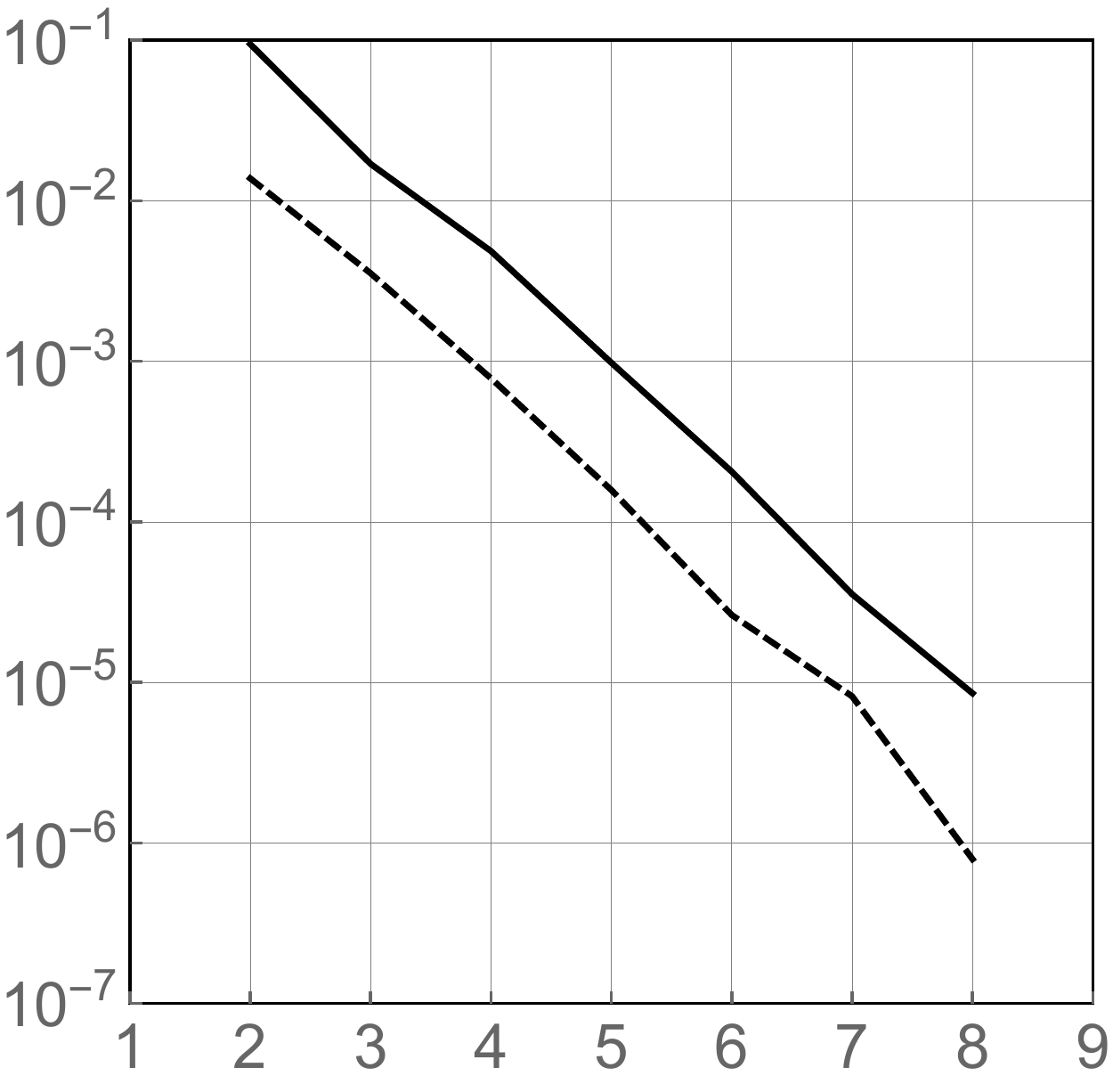}}\quad
\subfloat[Case 3.]{\includegraphics[width=0.30\textwidth]{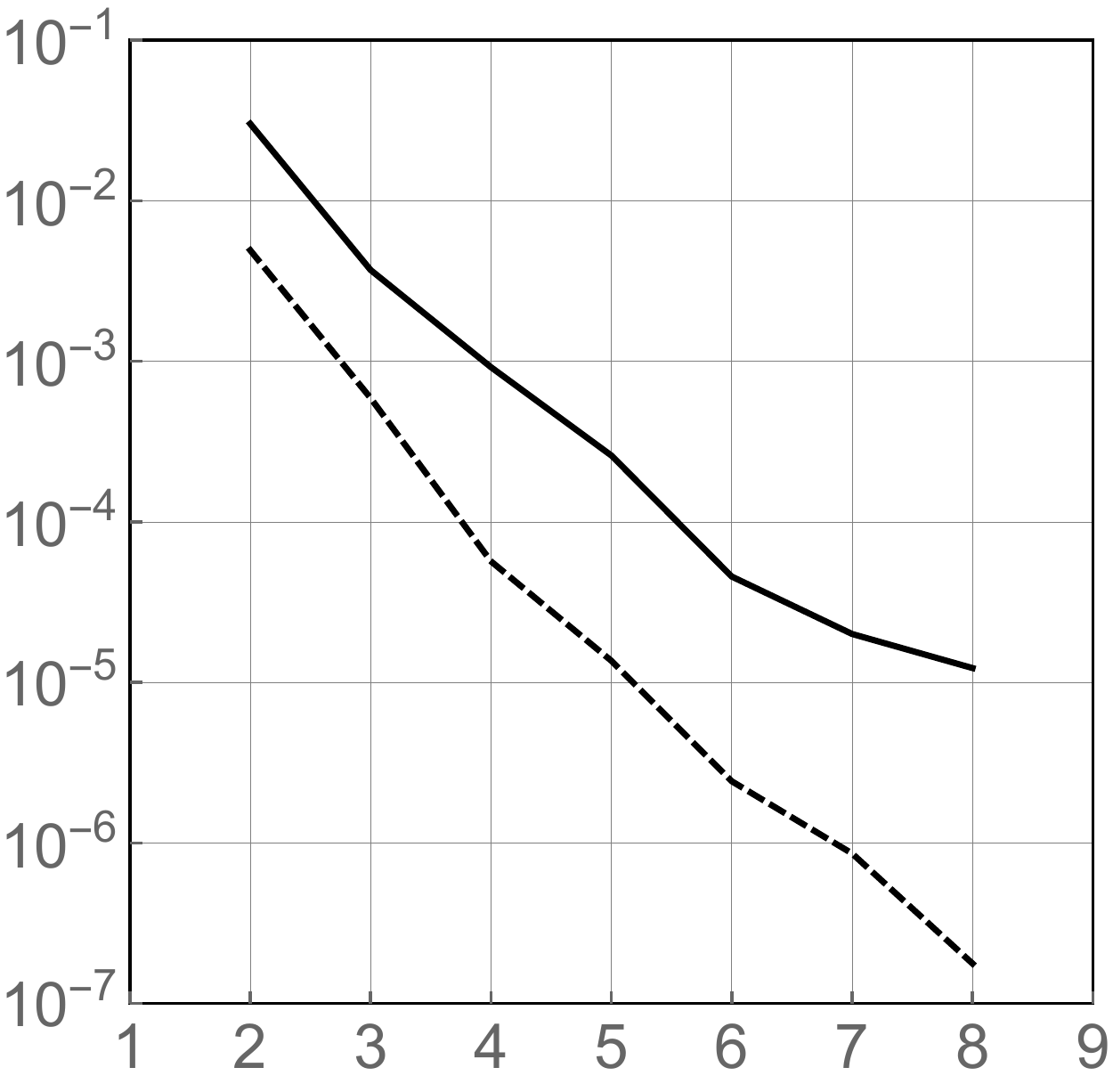}}\quad
\subfloat[Case 4.]{\includegraphics[width=0.30\textwidth]{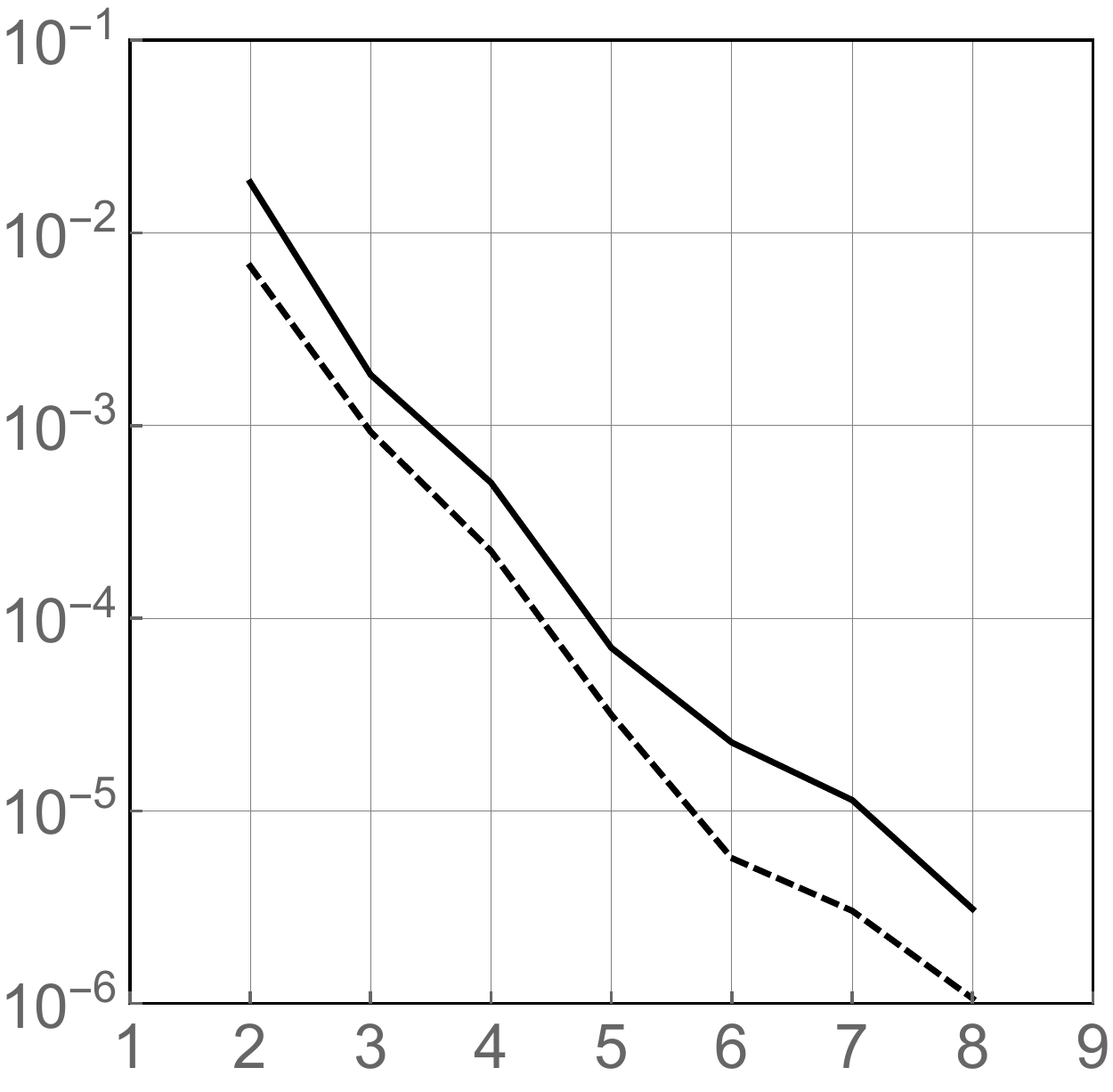}}
\caption{Quadrilaterals with Cusps: Cases 2--4: Estimated errors; $\log$-plot: Error vs $p$; Solid line = Reciprocal estimate, Dashed line = Auxiliary space estimate.}
\label{fig:cusperror}
\end{figure}

The numerical experiments differ from the previous cases since only the $p$-version is used.
In other words, the meshes of Figure~\ref{fig:cuspdomain} are used as is, without any $h$-refinement.
As is evident in the convergence and error estimation graphs of Figures~\ref{fig:cuspcase0}-\ref{fig:cusperror},
n the cases where the local angles close to $\pi/2$, exponential convergence is achieved, but
in the general case, when small geometric features are present, the convergence rates
stall to algebraic and not exponential.

\begin{table}[ht]
\caption{
   Tests on quadrilaterals with cusps.
}\label{tbl:eremenko}
\renewcommand\arraystretch{1}
\noindent

\begin{center}\footnotesize
\begin{tabular}{|l|l|l|l|l|l|}
\hline
Case & Method & Errors & Sizes & $M(Q_s)$\\
\hline
1 & $p$, $p=16$ & 9 & 1089 & 2.7823418091539533\\
2 & $p$, $p=16$ & 9 & 1633 & 1.8247899464782131\\
3 & $p$, $p=16$ & 7 & 2945 & 1.7864319361374579\\
4 & $p$, $p=16$ & 8 & 2945 & 0.8852475766134157\\
\hline
\end{tabular}
\
\end{center}
\end{table}
\section{Conclusions}

We have introduced a new class of ring domains, characterized by three parameters,
and given a formula for its modulus. By modifying the parameters, we obtain domains of increasing
computational challenge. For some specific sets of parameters we compute numerically
the modulus and compare the true error and two error estimates and show that the these two
error estimates behave in the same way as the true error. The computation is based on the
$hp$-FEM method and we show that nearly optimal convergence is obtained, when compared
to the theory of I. Babushka and his coauthors. 
This class of domains could be used for
benchmarking the numerical performance of FEM-software because of the scalability of computational challenge and 
the exactly known solution.

{\bf Acknowledgements.}
This research of the third author was supported by the Academy of Finland, Project 2600066611.

\end{document}